\documentclass[12pt]{article}
\usepackage{mathtools} 
\usepackage{amsmath, amssymb} 
\usepackage{amsthm}
\usepackage{enumerate}  % lets you do \begin{enumerate}[(a)] and such
	\usepackage{url} % formats url's nicely using \url{}
	\usepackage{authblk}
	\usepackage{tikz}
	\usepackage{tkz-graph}
	\usetikzlibrary{shapes}
	\usetikzlibrary{arrows}
	\usetikzlibrary{decorations.markings}
	
	\usepackage{graphicx}
	\usepackage{caption,subcaption}
	\usepackage{float}
	\usepackage{hyperref}
	
	\newcommand\cx{{\mathbb C}}% complexes

	\newcommand\ints{{\mathbb Z}}
	\newcommand\re{{\mathbb R}}%reals
	\newcommand\rats{{\mathbb Q}}

	\DeclarePairedDelimiter\abs{\lvert}{\rvert}%
	\DeclarePairedDelimiter\norm{\lVert}{\rVert}%
	
	\makeatletter
	\let\oldabs\abs
	\def\abs{\@ifstar{\oldabs}{\oldabs*}}
	\let\oldnorm\norm
	\def\norm{\@ifstar{\oldnorm}{\oldnorm*}}
	\makeatother
	
	\newcommand\sbs{\subseteq}

	\newcommand\comp[1]{{\mkern2mu\overline{\mkern-2mu#1}}}
	\newcommand\pmat[1]{\begin{pmatrix} #1 \end{pmatrix}}
	\newcommand\seq[4]{#1_{#2},#1_{#3},\ldots,#1_{#4}}

	\newtheoremstyle{plainsl}%
	{\topsep}
	{\topsep}
	{\slshape} % only non-default setting
	{}
	{\normalfont\bfseries}
	{.}
	{ }
	{}
	
	\swapnumbers
	
	{\theoremstyle{plainsl}
		\newtheorem{theorem}{Theorem}[section]
		\newtheorem{lemma}[theorem]{Lemma}
		\newtheorem{corollary}[theorem]{Corollary}}
	{\theoremstyle{remark}
		}
	
	\renewcommand\proof{\noindent\textsl{Proof. }}
	\newcommand\sqr[2]{{\vbox{\hrule height.#2pt
				\hbox{\vrule width.#2pt height#1pt \kern#1pt
					\vrule width.#2pt}\hrule height.#2pt}}}
	\renewcommand\qed{%
		\ifmmode\eqno\sqr53
		\else\nolinebreak\ \hfill\sqr53\medbreak\fi}

	\DeclareMathOperator{\supp}{supp}
	\DeclareMathOperator{\rk}{rk}
	\DeclareMathOperator{\tr}{tr}
	\DeclareMathOperator{\col}{col}

	\DeclareMathOperator{\adj}{adj}
	\DeclareMathOperator{\lcm}{lcm}

	\DeclareMathOperator{\arcs}{arcs}
		\DeclareMathOperator{\cl}{cl}
	
	\newcommand\ip[2]{\langle#1,#2\rangle}
	\newcommand\one{{\bf1}}

	\newcommand{\projonto}[2]{\operatorname{proj}_{#1}\!\left(#2\right)}

	\newenvironment{widetabular}{\tabular}{\endtabular}

	\usepackage{blkarray}

	\title{Subspace State Transfer in Coined Quantum Walks}
	\author{Yichi Xu and Hanmeng Zhan}
	\affil{Computer Science Department\\ Worcester Polytechnic Institute, Worcester, MA, USA\\\texttt{\{yxu10,hzhan\}@wpi.edu}}
	
\begin{document}

		\tikzset{->-/.style={decoration={
					markings,
					mark=at position .5 with {\arrow{>}}},postaction={decorate}}}
		\tikzset{lab/.style={draw=none,  inner sep = 2pt}}
		\maketitle
		
		\begin{abstract}
			We study a transport phenomenon in certain coined quantum walks where a subspace of states localized at a vertex gets transferred to another vertex. We first develop  characterizations for perfect and pretty good subspace state transfer using the spectral properties of a Hermitian weighted digraph obtained from the underlying graph. We then provide a polynomial-time algorithm that tests whether pointwise perfect subspace state transfer occurs at an integer step, given that the subspace and coins are rational. Finally, we construct several infinite families of examples that admit pointwise perfect $d$-dimensional subspace state transfer where $d\ge 2$.
		\end{abstract}
		
		\section{Introduction}
		The development of quantum algorithms has led to active research on transport phenomena in discrete quantum walks, such as perfect state transfer \cite{Lovett2010,Kendon2011,Kurzynski2011,Barr2014,Yalcnkaya2014,Stefanak2016,Stefanak2017,Zhan2019,Kubota2021,Guo2024a,Chen2024a,Kubota2025,Bhakta2024a}, periodicity \cite{Yoshie2019,Ito2020,Yoshie2023,Chen2024,Bhakta2024}, pretty good state transfer \cite{Chan2023,Zhan2025a,Bhakta2025}, and peak state transfer \cite{Guo2024}. While most existing studies focused on transferring a single state from one vertex to another, recently, Skoupy and Stefanak  \cite{Skoupy2025} observed that a discrete quantum walk on a graph with loops may simultaneously transfer two orthogonal states -- one supported by the loop, and one supported by the other outgoing arcs -- between two marked vertices, enabling transfer of a two-dimensional subspace of states from one vertex to another.

		In this paper, we study subspace state transfer on simple, unweighted graphs. We will consider arc-reversal quantum walks with reflection coins, and, using the spectral relation between these coins and certain Hermitian weighted digraph, we develop characterizations for perfect and pretty good subspace state transfer. These extend the theory of state transfer via Grover coins \cite{Lovett2010,Kendon2011,Kurzynski2011,Barr2014,Zhan2019,Kubota2021,Kubota2025,Bhakta2024a} or weighted Grover coins \cite{Chan2023,Zhan2025a}, where each coin is a reflection about a one-dimensional subspace.	
		
		Switching to coins that reflect about  higher-dimensional subspaces brings several advantages. For example, as shown by Kubota and Segawa \cite{Kubota2025}, the only $4$-regular circulant graphs that admit perfect state transfer in Grover coined walks are $X(\ints_{2m}, \pm\{a, m-a\})$ where $m$ is odd or oddly even, and the minimum transfer time is at least $m$. On the other hand, by Theorem \ref{circulant} of this paper, assigning certain higher-dimensional reflection coins to the sender and receiver vertices enables perfect state transfer on more circulant graphs ($X(\ints_{2m},\pm\{ a, m-a\})$ for any $m$), increases the dimension of subspaces of states transferred, and reduces the transfer time to a constant.
		
		\begin{center}		
			\begin{minipage}{0.5\textwidth}
				\centering
				\begin{tikzpicture}[
					scale=0.4,
					vertex/.style={circle,very thick,draw},
					lab/.style={draw=none, inner sep=1pt}
					]
					% vertices
					\node[vertex] (0) at (3.2,0) {};
					\node[vertex] (1) at (1.2,2.5) {};
					\node[vertex] (2) at (-1.2,2.5) {};
					\node[vertex,fill=blue] (3) at (-3.2,0) {};
					\node[vertex] (4) at (-1.2,-2.5) {};
					\node[vertex] (5) at (1.2,-2.5) {};
					
					% black edges
					\foreach \a/\b in {0/1,0/2,0/4,0/5,1/2,1/3,1/5,2/3,2/4,3/4,3/5,4/5}
					\draw[thick] (\a) -- (\b);
					
					% red directed edges with shifted labels
					\draw[->-, line width=0.5mm, blue]
					(3) -- node[lab, below, pos=0.35] {$i$} (1);
					
					% these two were overlapping the vertical black edge, so shift left
					\draw[->-, line width=0.5mm, blue]
					(3) -- node[lab, left,  pos=0.35, xshift=-3pt] {$1$} (2);
					
					\draw[->-, line width=0.5mm, blue]
					(3) -- node[lab, left,  pos=0.35, xshift=-3pt] {$-i$} (4);
					
					\draw[->-, line width=0.5mm, blue]
					(3) -- node[lab, above, pos=0.35, xshift=-3pt] {$-1$} (5);
				\end{tikzpicture}

				\captionof{figure}{At time $0$}
			\end{minipage}%
			\begin{minipage}{0.5\textwidth}
				\centering
				\begin{tikzpicture}[
					scale=0.4,
					vertex/.style={circle,very thick,draw},
					lab/.style={draw=none, inner sep=1pt}
					]
					% vertices
					\node[vertex,fill=red] (0) at (3.2,0) {};
					\node[vertex] (1) at (1.2,2.5) {};
					\node[vertex] (2) at (-1.2,2.5) {};
					\node[vertex] (3) at (-3.2,0) {};
					\node[vertex] (4) at (-1.2,-2.5) {};
					\node[vertex] (5) at (1.2,-2.5) {};
					
					% black edges
					\foreach \a/\b in {0/1,0/2,0/4,0/5,1/2,1/3,1/5,2/3,2/4,3/4,3/5,4/5}
					\draw[thick] (\a) -- (\b);
					
					% red directed edges: mirror of your left figure
					% (left had 3 -> 1 : i)  --> now 0 -> 2 : i
					\draw[->-, line width=0.5mm, red]
					(0) -- node[lab, below, pos=0.35] {$-i$} (2);
					
					% (left had 3 -> 2 : 1)  --> now 0 -> 1 : 1
					\draw[->-, line width=0.5mm, red]
					(0) -- node[lab, right, pos=0.35, xshift=3pt] {$-1$} (1);
					
					% (left had 3 -> 4 : -i) --> now 0 -> 5 : -i
					\draw[->-, line width=0.5mm, red]
					(0) -- node[lab, right, pos=0.35, xshift=3pt] {$i$} (5);
					
					% (left had 3 -> 5 : -1) --> now 0 -> 4 : -1
					\draw[->-, line width=0.5mm, red]
					(0) -- node[lab, above, pos=0.35, xshift=3pt] {$1$} (4);
					
				\end{tikzpicture}
				\captionof{figure}{At time $4$}
			\end{minipage}
		\end{center}
		
			\begin{center}		
			\begin{minipage}{0.5\textwidth}
				\centering
				\begin{tikzpicture}[
					scale=0.4,
					vertex/.style={circle,very thick,draw},
					lab/.style={draw=none, inner sep=1pt}
					]
					% vertices
					\node[vertex] (0) at (3.2,0) {};
					\node[vertex] (1) at (1.2,2.5) {};
					\node[vertex] (2) at (-1.2,2.5) {};
					\node[vertex,fill=blue] (3) at (-3.2,0) {};
					\node[vertex] (4) at (-1.2,-2.5) {};
					\node[vertex] (5) at (1.2,-2.5) {};
					
					% black edges
					\foreach \a/\b in {0/1,0/2,0/4,0/5,1/2,1/3,1/5,2/3,2/4,3/4,3/5,4/5}
					\draw[thick] (\a) -- (\b);
					
					% red directed edges with shifted labels
					\draw[->-, line width=0.5mm, blue]
					(3) -- node[lab, below, pos=0.35, xshift=-3pt] {$-i$} (1);
					
					% these two were overlapping the vertical black edge, so shift left
					\draw[->-, line width=0.5mm, blue]
					(3) -- node[lab, left,  pos=0.35, xshift=-3pt] {$1$} (2);
					
					\draw[->-, line width=0.5mm, blue]
					(3) -- node[lab, left,  pos=0.35, xshift=-3pt] {$i$} (4);
					
					\draw[->-, line width=0.5mm, blue]
					(3) -- node[lab, above, pos=0.35, xshift=-3pt] {$-1$} (5);
				\end{tikzpicture}

				\captionof{figure}{At time $0$}
			\end{minipage}%
			\begin{minipage}{0.5\textwidth}
				\centering
				\begin{tikzpicture}[
					scale=0.4,
					vertex/.style={circle,very thick,draw},
					lab/.style={draw=none, inner sep=1pt}
					]
					% vertices
					\node[vertex,fill=red] (0) at (3.2,0) {};
					\node[vertex] (1) at (1.2,2.5) {};
					\node[vertex] (2) at (-1.2,2.5) {};
					\node[vertex] (3) at (-3.2,0) {};
					\node[vertex] (4) at (-1.2,-2.5) {};
					\node[vertex] (5) at (1.2,-2.5) {};
					
					% black edges
					\foreach \a/\b in {0/1,0/2,0/4,0/5,1/2,1/3,1/5,2/3,2/4,3/4,3/5,4/5}
					\draw[thick] (\a) -- (\b);
					
					% red directed edges: mirror of your left figure
					% (left had 3 -> 1 : i)  --> now 0 -> 2 : i
					\draw[->-, line width=0.5mm, red]
					(0) -- node[lab, below, pos=0.35] {$i$} (2);
					
					% (left had 3 -> 2 : 1)  --> now 0 -> 1 : 1
					\draw[->-, line width=0.5mm, red]
					(0) -- node[lab, right, pos=0.35, xshift=3pt] {$-1$} (1);
					
					% (left had 3 -> 4 : -i) --> now 0 -> 5 : -i
					\draw[->-, line width=0.5mm, red]
					(0) -- node[lab, right, pos=0.35, xshift=3pt] {$-i$} (5);
					
					% (left had 3 -> 5 : -1) --> now 0 -> 4 : -1
					\draw[->-, line width=0.5mm, red]
					(0) -- node[lab, above, pos=0.35, xshift=1pt] {$1$} (4);
					
				\end{tikzpicture}
				\captionof{figure}{At time $4$}
			\end{minipage}
		\end{center}

		Our characterization connects subspace state transfer to properties of certain Hermitian weighted digraph. This leads to a polynomial-time algorithm that tests whether pointwise perfect subspace state transfer occurs at an integer step, given that the subspace and coins are rational. By marking vertices in Grover coined walks, we construct several infinite families of examples that admit pointwise perfect subspace state transfer.
		
		\section{Reflection coins and weighted digraphs} \label{spectral}
		
		Let $X$ be a connected graph on $n$ vertices. A discrete quantum walk on $X$ takes place on the \textsl{arcs} of $X$:
		\[\{(u,v): u,v \in V(X),\; u\sim v\}.\]
		We will refer to $u$ and $v$ as the \textsl{tail} and \textsl{head} of $(u,v)$, respectively. A \textsl{quantum state} associated with $X$ is a function from $\arcs(X)$ to $\cx$, usually represented by a unit vector in $\cx^{\arcs(X)}$. The transition matrix $U$ of the quantum walk is a product of two unitary matrices: the \textsl{arc-reversal operator} $R$, which sends each arc $(u,v)$ to its reverse $(v,u)$, and a \textsl{coin operator}
		\[C=\pmat{C_1 & & & \\ & C_2 & & \\ & & \ddots & \\ & & & C_n},\]
		where each $C_u$, called the \textsl{coin} at vertex $u$, is a $\deg(u)\times \deg(u)$ unitary matrix that acts on complex functions on the outgoing arcs of $u$. Implicit in $C_u$ is a linear order on the neighbors of $u$:
		\[\sigma_u: \{1,2,\cdots \deg(u)\} \to \{v: u\sim v\}.\]
		Given the state $e_{(u, \sigma_u(j))}$, the transition matrix first sends it to a linear combination of all arcs leaving $u$:
		\[C e_{(u, \sigma_u(j))} = (C_u) e_j = \sum_{\ell=1}^{\deg(u)} (C_u)_{\ell, j} e_{(u, \sigma_u(\ell))},\]
		and then swaps the coefficients on opposite arcs:
		\[R C e_{(u, \sigma_u(j))} = \sum_{\ell=1}^{\deg(u)} (C_u)_{\ell, j} e_{(\sigma_u(\ell), u)}.\]
		For more discussion on this model, see \cite[Section 2.1]{Godsil2017}.

\begin{figure}[h]
	\centering
	\resizebox{0.8\hsize}{!}{	\begin{widetabular}{c|c}
			coin operator & Hermitian weighted digraph\\
			\hline
			$\footnotesize \pmat{
				\frac{2}{3}J-I & & \\
				& \frac{2}{3}J-I &\\
				&&\ddots }$
			& 
			\raisebox{-.4\totalheight}{\includegraphics[width=2.2cm]{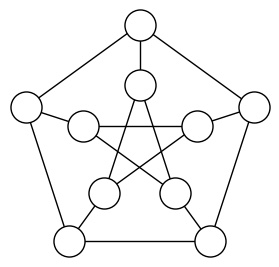}}
			\\
			\hline
			$\footnotesize\pmat{
				\color{blue}{-I} & & \\
				& \frac{2}{3}J-I &\\
				&&\ddots }$
			& \raisebox{-.4\totalheight}{\includegraphics[width=2.2cm]{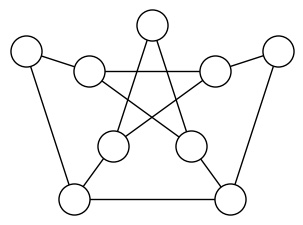}}\\
			\hline
			$\footnotesize \pmat{
				\color{blue}{\tiny \frac{2}{3}\pmat{1& i & i\\
						-i & 1 & 1 \\
						-i & 1 & 1}-I} & &\\
				& \frac{2}{3}J-I &\\
				& & \ddots }$
			&
			\raisebox{-.4\totalheight}{\includegraphics[width=2.2cm]{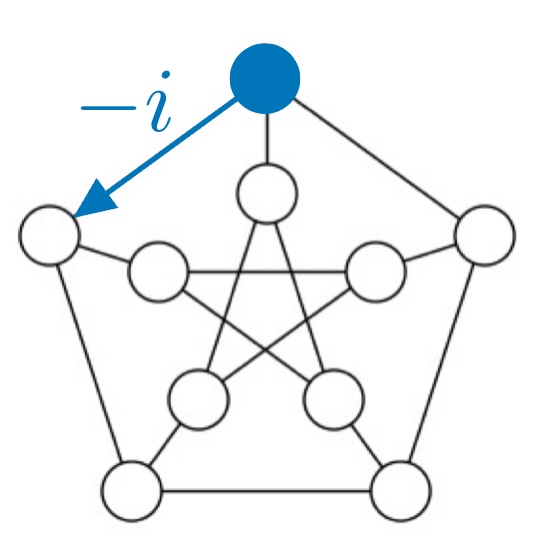}}\\
			\hline 
			$\footnotesize \pmat{
				\color{blue}{I-\frac{2}{3}J} & &  \\
				& \frac{2}{3}J-I &\\
				&&\ddots }$
			& \raisebox{-.4\totalheight}{\includegraphics[width=2.2cm]{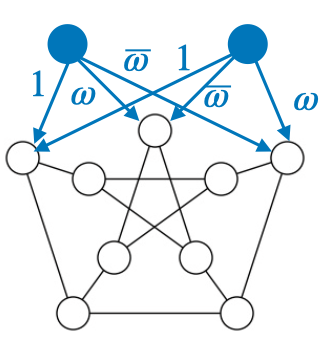}}
	\end{widetabular}}
	\caption{Reflection coins and Hermitian weighted digraphs; $\omega=e^{2\pi i/3}$}
	\label{reflection}
\end{figure}

		Throughout, we will assume each coin $C_u$ is a reflection:
		\[C_u^* = C_u,\quad C_u^2=I.\]
		Thus we may write 
		\[C_u = 2N_u N_u^*-I\]
		for some matrix $N_u$ with orthonormal columns. The rows of $N_u$ are indexed by the outgoing arcs of $u$, and the columns of $N_u$, if any, can be seen as indexed by some \textsl{clones} of $u$. We can then write 
		\[C=2NN^*-I,\]
		where
		\[N = \pmat{N_1 & & & \\ & N_2 & & \\ & & \ddots & \\ & & & N_n}\]
		is some weighted incidence matrix of the arcs of $X$ and the clones of vertices in $X$. The spectrum of $U=RC$ is closely related to that of the matrix $H=N^* R N$. Here $H$ also carries some combinatorial meaning: we can think of it as the Hermitian adjacency matrix of some digraph $Y$ constructed from $X$, where the vertices of $Y$ are the clones of $X$, and the weights on the arcs between these clones are determined by the entries of $N$.  Figure \ref{reflection} illustrates some reflection coin operators on the Peterson graph and their associated Hermitian weighted digraphs, with all weights scaled by $3$. We will describe these digraphs in more detail in Section \ref{sst}.

		The following spectral relation between $U$ and $H$ is a consequence of \cite[Chapter 2]{Godsil2023}. We omit the proof as it is analogous to that of \cite[Theorem 3.3]{Chan2023}.

		\begin{theorem}\label{thm:eprojs}\cite[Chapter 2]{Godsil2023}
			Let $U=R(2NN^*-I)$ be the transition matrix where $N$ has orthonormal columns. Let $H=N^*RN$. The eigenvalues of $U$ are $1$, $-1$ and $e^{\pm i\theta}$, where $\theta=\arccos\lambda$ for some eigenvalue $\lambda\in(-1,1)$ of $H$.  Moreover, if $F_{\theta}$ denotes the $e^{i\theta}$-eigenprojection of $U$, and $E_{\lambda}$ denotes the $\lambda$-eigenprojection of $H$, then $F_{\theta}$ and $E_{\lambda}$ are related in the following way.
			\begin{enumerate}[(i)]
				\item If $\theta=0$, then $\lambda=1$ and
				\[ N F_0 N^* = E_1.\]
				\item If $\theta=\pi$, then $\lambda=-1$ and
				\[N F_{\pi}N^* = E_{-1}.\]
				\item If $\theta \in(-\pi,0)\cup(0,\pi)$, then
				\[N F_{\theta}N^* = \frac{1}{2}E_{\lambda}.\]
			\end{enumerate}
		\end{theorem}

		\section{Subspace state transfer}\label{sst}
		
		Consider a quantum walk on a graph $X$ with transition matrix $U=RC$, where $C$ is a reflection. Let $a$ be a vertex of $X$. We say a state $x$ is a \textsl{coin state} at $a$ if 
		\[x = \sum_{u\sim a} w_u e_{(a,u)}\]
		for some $w$ satisfying $C_aw =w$; such state will be denoted $x_a(w)$. A basis for $\col(C+I)$ consisting of coin states is called a \textsl{coin basis}. Every coin basis $\mathcal{B}$ defines a Hermitian weighted digraph: if $M$ is the matrix with $\mathcal{B}$ as columns, then the Hermitian matrix $M^*RM$ encodes the adjacency matrix of some digraph $Y$ obtained from $X$ by cloning vertices and reweighting arcs based on $\mathcal{B}$. More specifically, each vertex $a$ of $X$ is cloned $\rk(C_a+I)$ times, and for every edge $\{a,b\}$ of $X$, the clones of $a$ and $b$ in $Y$ induce a weighted complete bipartite digraph, with weight 
		\[\ip{RM_b e_{\ell}}{M_a e_j}\]
		on the arc from the $j$-th clone of $a$ to the $\ell$-th clone of $b$.
		We will call $Y$ the \textsl{digraph relative to $\mathcal{B}$}, and write $H=M^*RM$.

		We are interested in transferring a subspace of coin states from one vertex to another. Given two vertices $a$ and $b$, we say there is \textsl{perfect subspace state transfer} from $a$ to $b$ if there exist a time $t \in \re$ and  a subspace $W$ of $\col(C_a+I)$ such that for any vector $w\in W$, the state $U^t x_a(w)$ is a coin state at $b$. Alternatively, let $x_a(W)$ denote the subspace
		\[x_a(W) = \mathrm{span}\{x_a(w): w\in W\}.\]
		There is perfect subspace state transfer from $a$ to $b$ if there are subspaces $W\sbs \col(C_a+I)$ and $V\sbs \col(C_b+I)$ such that $U^t x_a(W) = x_b(V)$ for some time $t$.  Note that this implies $\dim(W) = \dim(V)$. When $W$ and $V$ are known, we call this \textsl{perfect $(W,V)$-transfer from $a$ to $b$}.  When only the dimension $d$ of these subspaces is known, we call this \textsl{perfect $d$-dimensional subspace state transfer from $a$ to $b$}.		

		We will also consider a relaxation of this phenomenon called \textsl{pretty good subspace state transfer}, that is, there exists a subspace $W$ of $\col(C_a+I)$ such that for any $\epsilon>0$ and any vector $w\in W$, there is a time $t$ such that $U^t x_a(w)$ is $\epsilon$-close to some coin state at $b$:
		\[\norm{U^t x_a(w) - x_b(v)}<\epsilon,\quad v\in \col(C_b+I).\]

		Any basis for $x_a(W)\cup x_b(V)$	extends to a coin basis, and relative to this basis, $\dim(W)$ clones of $a$ arise from $W$, and $\dim(V)$ clones of $b$ arise from $V$. For convenience, we will refer to such basis as a \textsl{coin basis induced by $(W,V)$}, and those special clones as $W$-clones of $a$ and $V$-clones of $b$, respectively. We develop a characterization for perfect subspace state transfer in terms of these clones. For any subset $S$ of vertices in $Y$, let $B_S$ be the matrix whose columns are the standard basis vectors corresponding to the vertices in $S$. For any two matrices $M$ and $N$ with the same size, let $\ip{M}{N}$ denote their matrix inner product, that is,
		\[\ip{M}{N} = \tr(N^*M).\]
		The following lemma reduces the size of matrices involved in perfect subspace state transfer.

		\begin{lemma}\label{proj_tran}
			Let $W\sbs \col(C_a+I)$ and $V\sbs \col(C_b+I)$ be two subspaces with equal dimension. The following statements are equivalent.
			\begin{enumerate}[(i)]
				\item There is perfect $(W,V)$-transfer from $a$ to $b$ at time $t$.
				\item There is a matrix $N$ whose columns form an orthonormal coin basis induced by $(W,V)$ such that
						\[U^t NB_S = NB_T,\]
				where $S$ is the set of $W$-clones and $T$ is the set of $V$-clones.
				\item There is a matrix $N$ whose columns form an orthonormal coin basis induced by $(W,V)$ such that
				\[N^*U^t NB_S = B_T,\]
				where $S$ is the set of $W$-clones and $T$ is the set of $V$-clones.
			\end{enumerate}
		\end{lemma}
		\proof
		Suppose first that (i) holds. Let $d=\dim(W)$ and fix  an orthonormal basis  $\{\seq{w}{1}{2}{d}\}$ for $W$. Then
		\[U^t \pmat{x_a (w_1) & x_a(w_2) & \cdots & x_a(w_d)} = \pmat{x_b(v_1) & x_b(v_2) & \cdots & x_b(v_d)}\]
		for some orthonormal basis $\{\seq{v}{1}{2}{d}\}$ for $V$. Extend 		
		\[\{\seq{w}{1}{2}{d},\seq{v}{1}{2}{d}\}\]
		to an orthonormal coin basis, and let $N$ be the matrix with this basis as columns. Then 
		\begin{equation}\label{U^tN}
			U^t NB_S = NB_T.
		\end{equation}
		This proves (ii).
		
		Next suppose (ii) holds. Multiplying both sides of 
		\[U^t N B_S = NB_T,\]
		by $N^*$ on the left gives
\[
			N^*U^t NB_S = N^*NB_T=B_T.
\]
		This proves (iii).
		
		Finally suppose (iii) holds. Then
		\begin{equation}\label{reduction}
			B_T N^* U^t N B_S = B_T.
		\end{equation}
		Taking the trace on both sides yields
		\[\ip{U^t N B_S}{NB_T} = \ip{B_T}{B_T}.\]
		Thus the following Cauchy-Schwarz inequality is tight:
		\[\abs{\ip{U^t N B_S}{NB_T}}^2 \le \ip{B_S}{B_S}\ip{B_T}{B_T} = \ip{B_T}{B_T}^2\]
		and $U^t N B_S $ must be a scalar multiple of $N B_T$. It follows that for any $w\in W$,
		\[U^t x_a(w) \in \col(U^t NB_S) \sbs \col(NB_T) = x_a(V).\]
		 This proves (i).
		\qed

		Let $f_t$ be the function defined on the interval $[-1,1]$ by
		\[f_t(x) = \cos(t\arccos(x)).\]
		Note that when $t$ is a positive integer, $f_t(x)$ coincides with the Chebyshev polynomial of the first kind. For any Hermitian matrix $H$ with eigenvalues in $[-1,1]$ and spectral decomposition
		\[H = \sum_{\lambda} \lambda E_{\lambda},\]
 we define 
		\[f_t(H) = \sum_{\lambda} f_t(\lambda) E_{\lambda}.\]
		In particular, if $H$ is the Hermitian adjacency matrix relative to some orthonormal coin basis:
		\[H = N^*RN,\]
		then we can expand $N^*U^tN$ using the eigenvalues $e^{i\theta_r}$ and eigenprojections $F_r$ of $U$:
		\[N^* U^t N = \sum_r e^{i t\theta_r} N^*F_r N .\]
 	By Theorem \ref{thm:eprojs}, the right hand side reduces to
		\[\sum_{\lambda} \cos(t\arccos\lambda) E_{\lambda} = f_t(H).\]
		Thus Lemma  \ref{proj_tran} (iii) can be restated in terms of $f_t(H)$. We remark that similar approaches to one-dimensional subspace state transfer have appeared in several works \cite{Zhan2019,Kubota2021,Chan2023,Bhakta2024a,Kubota2025,Guo2024,Zhan2025a,Bhakta2025}.
		
		\begin{lemma}
			Let $W\sbs \col(C_a+I)$ and $V\sbs \col(C_b+I)$ be two subspaces with dimension $d$. There is perfect $(W,V)$-transfer from $a$ to $b$ at time $t$ if and only if 
			\[f_t(H) B_S = B_T,\]
			or equivalently,
			\[\ip{f_t(H) B_S}{B_T} = d,\]
			relative to some orthonormal coin basis induced by $(W,V)$. 
		\end{lemma}
		
		We also obtain an equivalent formulation for pretty good subspace state transfer in terms of $f_t(H)$.
		
		\begin{lemma}\label{pgsst}
			Let $W\sbs \col(C_a+I)$ and $V\sbs \col(C_b+I)$ be two subspaces with equal dimension. There is pretty good $(W,V)$-transfer from $a$ to $b$ at time $t$ if and only if for any $\epsilon>0$, there is a time $t$ such that
			\[\norm{f_t(H) B_S - B_T}<\epsilon,\]
			or equivalently,
			\[\abs{\ip{f_t(H) B_S}{B_T} -d}<\epsilon,\]
			relative to some orthonormal coin basis induced by $(W,V)$. 
		\end{lemma}

\section{A characterization}\label{sec:char}
	
	To find necessary and sufficient conditions for subspace state transfer, we consider the inner product $\ip{f_t(H)B_S}{B_T}$ with $\abs{S}=\abs{T}$. We have
	\begin{align}
		\abs{\ip{f_t(H)B_S}{B_T}} & = \abs{\sum_{\lambda} f_t(\lambda) \ip{E_{\lambda} B_S}{B_T}}\notag\\
		& \le \sum_{\lambda} \abs{f_t(\lambda)} \abs{\ip{E_{\lambda} B_S}{B_T}} \label{ineq:first}\\
		&\le\sum_{\lambda}\abs{\ip{E_{\lambda} B_S}{B_T}}\label{ineq:second}\\
		&=\sum_{\lambda} \abs{\ip{E_{\lambda} B_S}{E_{\lambda}B_T}} \notag \\
		&\le \sum_{\lambda}\sqrt{\ip{E_{\lambda} B_S}{E_{\lambda}B_S}\ip{E_{\lambda} B_T}{E_{\lambda}B_T}}\label{ineq:third}\\
		&\le \sqrt{\sum_{\lambda}\ip{E_{\lambda} B_S}{E_{\lambda}B_S}}\sqrt{\sum_{\lambda}\ip{E_{\lambda} B_T}{E_{\lambda}B_T}}\label{ineq:fourth}\\
		&=\sqrt{\sum_{\lambda}\ip{B_S}{B_S}}\sqrt{\sum_{\lambda}\ip{ B_T}{B_T}} \notag \\
		&=\sqrt{\abs{S}}\sqrt{\abs{T}}\notag\\
		&=\ip{B_T}{B_T}. \notag 
		\notag        \end{align}		
	Note that Inequaltiy \eqref{ineq:first} is tight if and only if all the nonzero products $f_t(\lambda)\ip{E_{\lambda} B_S}{B_T}$ have the same complex argument. Inequality \eqref{ineq:second} is tight if and only if $f_t(\lambda)=\pm 1$. Inequality \eqref{ineq:third} is tight if and only if $E_{\lambda} B_S$ is a scalar multiple of $E_{\lambda} B_T$. Inequality \eqref{ineq:fourth} is tight if and only if $\ip{E_{\lambda}B_S}{B_S} = \ip{E_{\lambda}B_T}{B_T}$. These motivate the following definitions.

		Let $Y$ be a Hermitian weighted digraph, with spectral decomposition
		\[H = \sum_{\lambda} \lambda E_{\lambda}.\]
		The \textsl{eigenvalue support} of a subset $S$ of vertices in $Y$ is the set
		\[\Lambda_S = \{\lambda: E_{\lambda} B_S\ne 0\}.\]
		We say two subsets $S$ and $T$ of vertices in $Y$ are \textsl{strongly-cospectral} if for each eigenvalue $\lambda$ of $H$,
		\[E_{\lambda} B_S = \pm  E_{\lambda} B_T.\]
		If $S$ and $T$ are strongly-cospectral, then $\Lambda_S=\Lambda_T$, which can be partitioned into two disjoint sets
		\[\Lambda_{S,T}^+ = \{\lambda: E_{\lambda} B_S =  E_{\lambda} B_T\ne 0\}\]
		and 
		\[\Lambda_{S,T}^- = \{\lambda: E_{\lambda} B_S = - E_{\lambda} B_T\ne 0\}.\]

		We now give a characterization of perfect subspace state transfer in terms of the eigenvalues and eigenprojections of $H$.
		
		\begin{theorem}\label{psst_char}
			Let $W\sbs \col(C_a+I)$ and $V\sbs \col(C_b+I)$ be two subspaces with equal dimension. There is perfect $(W,V)$-transfer from $a$ to $b$ at time $t$ if and only if relative to some orthonormal coin basis induced by $(W,V)$, the Hermitian adjacency matrix $H$ satisfies the following.
			\begin{enumerate}[(i)]
				\item $S$ and $T$ are strongly-cospectral.
				\item  $t\arccos\lambda$ is an even multiple of $\pi$ for each eigenvalue $\lambda\in \Lambda_{S,T}^+$.
				\item $t\arccos\lambda$ is an odd multiple of $\pi$ for each eigenvalue $\lambda\in \Lambda_{S,T}^-$.
			\end{enumerate}
		\end{theorem}
		\proof
		Suppose first that $f_t(H)B_S = B_T$. Then $\ip{f_t(H)B_S}{B_T} = \ip{B_T}{B_T}$. Thus $\ip{f_t(H)B_S}{B_T}>0$, and all inequalities in \eqref{ineq:first}, \eqref{ineq:second}, \eqref{ineq:third}, \eqref{ineq:fourth} are tight. In particular, tightness of \eqref{ineq:first} now implies 
		\begin{equation}\label{tight3strong}
			f_t(\lambda)\ip{E_{\lambda} B_S}{B_T}\ge 0,\quad \forall \lambda.
		\end{equation}
		Let $\lambda \in \Lambda_S$. Since equality holds in \eqref{ineq:third}, there is a scalar $\alpha_{\lambda}$ such that
		\[E_{\lambda} B_S = \alpha_{\lambda} E_{\lambda} B_T.\] Plugging this into \eqref{tight3strong} shows 
		\[\alpha_{\lambda}f_t(\lambda) >0.\]
		Since \eqref{ineq:fourth} is tight, we must have $\alpha_{\lambda}=\pm 1$. Finally, by tightness of \eqref{ineq:second}, if $\alpha_{\lambda}=1$, then $f_t(\lambda)=1$, in which case $t\arccos\lambda$ is an even multiple of $\pi$, and if $\alpha_{\lambda}=-1$, then $f_t(\lambda)=-1$, in which case $t\arccos \lambda$ is an odd multiple $\pi$. This proves conditions (i) -- (iii).
		
		Conversely, suppose conditions (i) -- (iii) hold. Then
		\begin{align*}
			f_t(H)B_S &= \sum_{\lambda} f_t(\lambda) E_{\lambda} B_S\\
			&= \sum_{\lambda\in \Lambda_{S,T}^+} E_{\lambda} B_S  - \sum_{\lambda \in \Lambda_{S,T}^-} E_{\lambda} B_S\\
			&= \sum_{\lambda\in \Lambda_{S,T}^+} E_{\lambda} B_T  + \sum_{\lambda \in \Lambda_{S,T}^-} E_{\lambda} B_T\\
			&= \sum_{\lambda} E_{\lambda} B_T \\
			&=B_T. \tag*{\sqr53}
		\end{align*}
		
		Note that the characterization in Theorem \ref{psst_char} is symmetric in $a$ and $b$. 
		\begin{lemma}
			If there is perfect $(W,V)$-transfer from $a$ to $b$ at time $t$, then there is perfect $(W,V)$-transfer from $b$ to $a$ at time $t$.
		\end{lemma}

		We say $U$ is \textsl{subspace periodic} at $a$ if there is a time $t\in \re$ and a subspace $W$ of $\col(C_a+I)$ such that $x_a(W)$ is $U^t$-invariant; when $W$ is known, $U$ is said to be \textsl{$W$-periodic at $a$}. If, in addition, 
		\[U^t x_a(w)  =  x_a(w)\]
		for any $w\in W$, then $U$ is said to be \textsl{pointwise subspace periodic} at $a$. Setting $a=b$ and $W=V$ in Theorem \ref{psst_char} gives a characterization of pointwise subspace periodicity.

		\begin{theorem}\label{char_psp}
		Let $W$ be a subspace of $\col(C_a+I)$ with an orthonormal basis $\{\seq{w}{1}{2}{d}\}$.  Let $H$ be the Hermitian adjacency matrix relative to any orthonormal coin basis containing
		\[\{x_a(w_1), x_a(w_2), \cdots, x_a(w_d)\}.\]
		Then $U$ is pointwise $W$-periodic at $a$ with period $t$ if and only if $t\arccos\lambda$ is an even multiple of $\pi$ for each $\lambda\in \Lambda_S$.
		\end{theorem}

		The next corollary is an immediate consequence of the symmetry of perfect subspace state transfer.
		
		\begin{corollary}
			 If $U$ admits perfect $(W,V)$-transfer from $a$ to $b$ at time $t$, then it is pointwise $W$-periodic at $a$ and pointwise $V$-periodic at $b$, both with period $2t$.
		\end{corollary}
		
		This implies a ``monogamy" property of perfect subspace state transfer.
		
		\begin{corollary}
			If there is perfect $(W,V)$-transfer from $a$ to $b$ and perfect $(W, V')$-transfer from $a$ to $b'$, then $b=b'$ and $V=V'$.
		\end{corollary}

		To characterize pretty good subspace state transfer, we need two versions of Kronecker's approximation theorem, one for approximation in real steps and one for approximation in integer steps.

		\begin{theorem}\cite[Ch 3]{Levitan1983} \label{Levitan1983}
			Given $\seq{\alpha}{1}{2}{n}\in \re$ and $\seq{\beta}{1}{2}{n}\in \re$, the following are equivalent.
			\begin{enumerate}[(i)]
				\item For any $\epsilon>0$, the system
				\[\abs{q\alpha_r - \beta_r}<\epsilon \pmod{2\pi}, \quad r=1,2,\cdots, n\]
				has a solution $q\in\re$.
				\item For any set $\{\seq{\ell}{1}{2}{n}\}$ of integers such that 
				\[\ell_1 \alpha_1 + \cdots + \ell_n \alpha_n =0,\]
				we have
				\[\ell_1 \beta_1 + \cdots + \ell_n \beta_n \equiv 0\pmod{2\pi}.\]
			\end{enumerate}
		\end{theorem}

		\begin{theorem}\cite{Gonek2016}\label{thm:Kronecker}
			Given $\seq{\alpha}{1}{2}{n}\in \re$ and $\seq{\beta}{1}{2}{n}\in \re$, the following are equivalent.
			\begin{enumerate}[(i)]
				\item For any $\epsilon>0$, the system
				\[\abs{q\alpha_r - \beta_r-p_r}<\epsilon, \quad r=1,2,\cdots, n\]
				has a solution $\{q, \seq{p}{1}{2}{n}\} \in \ints^{n+1}$.
				\item For any set $\{\seq{\ell}{1}{2}{n}\}$ of integers such that 
				\[\ell_1 \alpha_1 + \cdots + \ell_n \alpha_n \in \ints,\]
				we have
				\[\ell_1 \beta_1 + \cdots + \ell_n \beta_n \in \ints.\]
			\end{enumerate}
		\end{theorem}

The following characterization for pretty good subspace state transfer extends Theorem 3.6 in \cite{Zhan2025a} on pretty good state transfer via weighted Grover coins.
\begin{theorem}\label{pgsst_char}
	Let $W\sbs \col(C_a+I)$ and $V\sbs \col(C_b+I)$ be two subspaces with equal dimension. There is pretty good $(W,V)$-transfer from $a$ to $b$ if and only if relative to some orthonormal coin basis induced by $(W,V)$, the Hermitian adjacency matrix $H$ satisfies the following.
	\begin{enumerate}[(i)]
		\item $S$ and $T$ are strongly-cospectral.
		\item  For any set $\{\ell_{\lambda}:\lambda\in\Lambda_S\}$ of integers such that
		\begin{equation}\label{K1}
			\sum_{\lambda \in \Lambda_S} \ell_{\lambda} \arccos \lambda =0,
		\end{equation}
		we have
		\[\sum_{\lambda\in\Lambda^-_{S,T}}\ell_{\lambda} \equiv 0\pmod2.\]
	\end{enumerate}
	Moreover, $x_b(V)$ can be approximated in integer steps if and only if (i) and (ii) hold with Equation \eqref{K1} replaced by
	\[\sum_{\lambda \in \Lambda_S} \ell_{\lambda} \arccos \lambda \equiv 0 \pmod{2\pi}.\]
\end{theorem}		
\proof
By Lemma \ref{pgsst}, there is pretty good $(W,V)$-transfer from $a$ to $b$ if and only if 
\[\{\ip{f_t(H) B_S}{B_T}: t\in \re\}\]
has $d$ as a limit point. 
It follows from the inequalities in  \eqref{ineq:third}  and \eqref{ineq:fourth} that $S$ and $T$ need to be strong-cospectrality to allow pretty good $(W,V)$-transfer. Given (i), set $\sigma_{\lambda}=0$ if $\lambda\in \Lambda_{S,T}^+$, and $\sigma_{\lambda}=1$ if $\lambda\in \Lambda_{S,T}^-$. Pretty good $(W,V)$-transfer occurs from $a$ to $b$ if and only if for any $\epsilon>0$, there is a $t\in \re$ such that
\[	\abs{f_t(\lambda)(-1)^{\sigma_{\lambda}} - 1}= \abs{\cos(t\arccos\lambda)- \cos(\sigma_{\lambda}\pi)}<\epsilon,\quad \forall \lambda \in \Lambda_{S},\]
or equivalently, the system
\[ \abs{t \arccos\lambda- \sigma_{\lambda}\pi - 2k_{\lambda}\pi}<\epsilon, \quad \forall \lambda \in \Lambda_S\]
has a solution $\{t\}\cup \{k_{\lambda}: \lambda\in \Lambda_S\}$ where $t\in \re$ and $k_{\lambda}\in \ints$. By Theorem \ref{Levitan1983}, this happens if and only if for any for any set  $\{\ell_{\lambda}:\lambda\in\Lambda_a\}$ of integers such that
\[\sum_{\lambda \in \Lambda_a} \ell_{\lambda} \arccos \lambda =0,\]
we have
\[\sum_{\lambda\in\Lambda_{ab}}\ell_{\lambda} \sigma_{\lambda}=\sum_{\lambda\in\Lambda_{ab}^-}\ell_{\lambda} \equiv 0\pmod2.\]
The characterization for pretty good $(W,V)$-transfer from $a$ to $b$ in integer steps follows similarly from Theorem \ref{thm:Kronecker}.
\qed

\section{Pointwise subspace state transfer}

Sometimes we may want to preserve the amplitude on  each arc during state transfer. We say there is \textsl{pointwise perfect subspace state transfer} from $a$ to $b$ if there is some time $t\in \re$ and some subspace $W\sbs \col(C_a+I) \cap \col(C_b+I)$ such that for any vector $w\in W$, there is a unimodular $\gamma$ such that 
\[U^t x_a(w) = \gamma x_b(w).\]
Most studies on perfect and pretty good state transfer via weighted Grover coins fall into this category; the following is an example with $W=\mathrm{span}\{\one\}$ and $\gamma=1$. When $W$ is known, this is called \textsl{pointwise perfect $W$-transfer from $a$ to $b$}. 
It follows from linearity of $U$ that  if this happens, then the phase $\gamma$ does not depend on the choice of $w$.

\begin{center}		
	\begin{minipage}{0.5\textwidth}
		\centering
		\begin{tikzpicture}[
			scale=0.4,
			vertex/.style={circle,very thick,draw},
			lab/.style={draw=none, inner sep=1pt}
			]
			% vertices
			\node[vertex] (0) at (3.2,0) {};
			\node[vertex] (1) at (1.2,2.5) {};
			\node[vertex] (2) at (-1.2,2.5) {};
			\node[vertex,fill=blue] (3) at (-3.2,0) {};
			\node[vertex] (4) at (-1.2,-2.5) {};
			\node[vertex] (5) at (1.2,-2.5) {};
			
			% black edges
			\foreach \a/\b in {0/1,0/2,0/4,0/5,1/2,1/3,1/5,2/3,2/4,3/4,3/5,4/5}
			\draw[thick] (\a) -- (\b);
			
			% red directed edges with shifted labels
			\draw[->-, line width=0.5mm, blue]
			(3) -- node[lab, below, pos=0.35] {$1$} (1);
			
			% these two were overlapping the vertical black edge, so shift left
			\draw[->-, line width=0.5mm, blue]
			(3) -- node[lab, left,  pos=0.35, xshift=-3pt] {$1$} (2);
			
			\draw[->-, line width=0.5mm, blue]
			(3) -- node[lab, left,  pos=0.35, xshift=-3pt] {$1$} (4);
			
			\draw[->-, line width=0.5mm, blue]
			(3) -- node[lab, above, pos=0.35,] {$1$} (5);
		\end{tikzpicture}

		\captionof{figure}{At time $0$}
	\end{minipage}%
	\begin{minipage}{0.5\textwidth}
		\centering
		\begin{tikzpicture}[
			scale=0.4,
			vertex/.style={circle,very thick,draw},
			lab/.style={draw=none, inner sep=1pt}
			]
			% vertices
			\node[vertex,fill=red] (0) at (3.2,0) {};
			\node[vertex] (1) at (1.2,2.5) {};
			\node[vertex] (2) at (-1.2,2.5) {};
			\node[vertex] (3) at (-3.2,0) {};
			\node[vertex] (4) at (-1.2,-2.5) {};
			\node[vertex] (5) at (1.2,-2.5) {};
			
			% black edges
			\foreach \a/\b in {0/1,0/2,0/4,0/5,1/2,1/3,1/5,2/3,2/4,3/4,3/5,4/5}
			\draw[thick] (\a) -- (\b);
			
			% red directed edges: mirror of your left figure
			% (left had 3 -> 1 : i)  --> now 0 -> 2 : i
			\draw[->-, line width=0.5mm, red]
			(0) -- node[lab, below, pos=0.35] {$1$} (2);
			
			% (left had 3 -> 2 : 1)  --> now 0 -> 1 : 1
			\draw[->-, line width=0.5mm, red]
			(0) -- node[lab, right, pos=0.35, xshift=3pt] {$1$} (1);
			
			% (left had 3 -> 4 : -i) --> now 0 -> 5 : -i
			\draw[->-, line width=0.5mm, red]
			(0) -- node[lab, right, pos=0.35, xshift=3pt] {$1$} (5);
			
			% (left had 3 -> 5 : -1) --> now 0 -> 4 : -1
			\draw[->-, line width=0.5mm, red]
			(0) -- node[lab, above, pos=0.35] {$1$} (4);
			
		\end{tikzpicture}
		\captionof{figure}{At time $6$}
	\end{minipage}
\end{center}

As a strengthening of subspace state transfer, pointwise $W$-transfer from $a$ to $b$ imposes conditions on a special type of orthonormal coin bases, each of which contains
\[\{x_a(w_1), x_a(w_2), \cdots, x_a(w_d), x_b(w_1), x_b(w_2), \cdots, x_b(w_d)\},\]
for some orthonormal basis $\{\seq{w}{1}{2}{d}\}$ of $W$. We will call such a basis an \textsl{orthonormal coin basis induced by $W$}.

\begin{lemma}
	Let $W$ be a subspace of $\col(C_a+I)\cap \col(C_b+I)$. Pointwise perfect $W$-transfer from $a$ to $b$ occurs at time $t$ if and only if there is a unimodular $\gamma\in \cx$ such that
	\[f_t(H) B_S = \gamma B_T\]
	relative to some (or any) orthonormal coin basis induced by $W$.
\end{lemma}

Thus if $S$ and $T$ are the $W$-clones of $a$ and the $W$-clones of $b$, respectively, then a necessary condition for pointwsie $W$-transfer is that, for some unimodular $\gamma$,
\[E_{\lambda} B_S = \pm \gamma E_{\lambda} B_T\]
for each eigenprojection $E_{\lambda}$ of the special Hermitian weighted digraph. Following  Section \ref{sec:char}, we say $S$ is \textsl{$\gamma$-strongly-cospectral} to $T$ if the above holds, and define
\[\Lambda_{S,T}^{\gamma} = \{\lambda: E_{\lambda} B_S =  \gamma E_{\lambda} B_T\ne 0\}\]
and 
\[\Lambda_{S,T}^{-\gamma} = \{\lambda: E_{\lambda} B_S = -\gamma E_{\lambda} B_T\ne 0\}.\]
By relaxing the notion of strong-cospectrality, we obtain a characterization of pointwise $W$-transfer relative to an arbitrary orthonormal coin basis induced by $W$.

\begin{theorem}\label{char-ptwise-psst}
		Let $W$ be a subspace of $\col(C_a+I)\cap \col(C_b+I)$. Pointwise perfect $W$-transfer from $a$ to $b$ occurs at time $t$ if and only if there is a unimodular $\gamma\in \cx$ such that relative to some (or any) orthonormal coin basis induced by $W$, the following hold.
		\begin{enumerate}[(i)]
			\item $S$ is $\gamma$-strongly-cospectral to $T$.
			\item  $t\arccos\lambda$ is an even multiple of $\pi$ for each eigenvalue $\lambda\in \Lambda_{S,T}^{\gamma}$.
			\item $t\arccos\lambda$ is an odd multiple of $\pi$ for each eigenvalue $\lambda\in \Lambda_{S,T}^{-\gamma}$.
		\end{enumerate}
\end{theorem}

Similarly, we say there is \textsl{pointwise pretty good subspace state transfer} from $a$ to $b$ if there is a subspace $W$ of $\col(C_a+I)\cap \col(C_b+I)$ such that for any $\epsilon>0$ and any vector $w\in W$, there exist a time $t$ and a unimodular $\gamma$ such that 
\[\norm{U^t x_a(w) - \gamma x_b(w)}<\epsilon.\]

\begin{theorem}\label{char-ptwise-pgsst}
	Let $W$ be a subspace of $\col(C_a+I)\cap \col(C_b+I)$. Pointwise pretty good $W$-transfer occurs from $a$ to $b$ if and only if there is a unimodular $\gamma\in\cx$ such that relative to some (or any) orthonormal coin basis induced by $W$, the following hold. 
	\begin{enumerate}[(i)]
		\item $S$ is $\gamma$-strongly-cospectral to $T$.
		\item  For any set $\{\ell_{\lambda}:\lambda\in\Lambda_S\}$ of integers such that
		\begin{equation}\label{K2}
			\sum_{\lambda \in \Lambda_S} \ell_{\lambda} \arccos \lambda =0,
		\end{equation}
		we have
		\[\sum_{\lambda\in\Lambda^{-\gamma}_{S,T}}\ell_{\lambda} \equiv 0\pmod2.\] 
	\end{enumerate}
		Moreover, the target states can be approximated in integer steps if and only if (i) and (ii) hold with Equation \eqref{K2} replaced by
	\[\sum_{\lambda \in \Lambda_S} \ell_{\lambda} \arccos \lambda \equiv 0 \pmod{2\pi}.\]
\end{theorem}

\section{$\gamma$-strong-cospectrality}		\label{gamma-sc}
In this section, we provide a characterization for $\gamma$-strong-cospectrality using certain rational functions.  Much of the theory can be seen as an extension of \cite{Godsil2024} from two vertices to two sets of vertices. We will use these results to develop a polynomial-time algorithm that tests whether pointwise perfect subspace state transfer occurs at an integer step, given the subspace and coins are rational, in Section  \ref{ptwise_psst}. This complements the polynomial-time algorithm for perfect state transfer in continuous quantum walks \cite{Coutinho2017b}.

For any matrix $M$, subset $S$ of its row indices, and subset $T$ of its column indices, let $M_{S,T}$ denote the submatrix of $M$ with rows in $S$ and columns in $T$, and let $M\backslash[S, T]$ denote the matrix obtained from $M$ by deleting the rows in $S$ and the columns in $T$.

Let $Y$ be any Hermitian weighted digraph. Let $H$ be its Hermitian adjacency matrix with spectral decomposition
\[H = \sum_{\lambda} \lambda E_{\lambda}.\]
%Let $\psi(x)$ denote the characteristic polynomial of $H$, that is,
%\[	\psi(x) = \det(xI-H).\]
For any two subsets $S$ and $T$ of vertices in $Y$ with equal size, define
\[
	\psi_{S,T}(x) = \tr((xI-H)^{-1}_{S,T}).
\]
Note that $\psi_{S,T}(x)$ is a rational function with simple poles, as 
\begin{equation}\label{adjdet}
	\sum_{\lambda} \frac{1}{x-\lambda} (E_{\lambda})_{S,T} = (xI-H)^{-1}_{S,T}=\frac{\adj(xI-H)_{S,T}}{\det(xI-H)}.
\end{equation}
When $S=T$, we will abbreviate $\psi_{S,T}(x) $ as $\psi_S(x)$. Since $E_{\lambda}$ is positive semidefinite, $E_{\lambda}B_S=0$ if and only if $\tr((E_{\lambda})_{S,S})=0$. Thus the eigenvalue support of $S$ can be read off from the poles of $\psi_S(x)$. 

\begin{lemma}\label{evsp}
Let $S$ be any subset of vertices in a Hermitian weighted digraph. The eigenvalue support of $S$ consists of the poles of $\psi_S(x)$.
\end{lemma}

We say $S$ and $T$ are \textsl{cospectral} if 
\[\psi_S(x) = \psi_T(x).\] 
The following lemma offers two alternative definitions of cospectrality.

\begin{lemma}
Let $Y$ be a Hermitian weighted digraph with spectral decomposition
\[H = \sum_{\lambda} \lambda E_{\lambda}.\]
Let $S$ and $T$ be two subsets of vertices in $Y$ with equal size. The following are equivalent. 
	\begin{enumerate}[(i)]
		\item $S$ and $T$ are cospectral.
		\item For each $\lambda$,
		\[\tr((E_{\lambda})_{S,S}) = \tr((E_{\lambda})_{T,T}).\]
		\item For each nonnegative integer $m$,
		\[\tr(H^m_{S,S}) = \tr(H^m_{T,T}).\]
	\end{enumerate}
\end{lemma}
\proof
The equivalence between (i) and (ii) follows from Equation \eqref{adjdet}. The equivalence between (ii) and (iii) follows from the spectral decomposition of $H$, and the fact that each $E_{\lambda}$ is a polynomial in $H$.
\qed

The next result gives a characterization for $\gamma$-strong-cospectrality.

\begin{theorem}\label{gamma-sc-char}
	Let $Y$ be a Hermitian weighted digraph. Let $S$ and $T$ be two subsets of vertices in $Y$ of equal size. Then $S$ is $\gamma$-strongly-cospectral to $T$ for some unimodular complex number $\gamma$ if and only if the following hold.
	\begin{enumerate}[(i)]
		\item $S$ and $T$ are cospectral.
		\item $\gamma\psi_{S,T}(x)$ is a rational function over $\re$.
		\item The set of poles of $\psi_S(x)$ is the union of 
		\begin{equation}\label{Gamma+}
			\{\text{poles of $\psi_S(x)$ but not of $\psi_S(x) - \gamma \psi_{S,T}(x)$}\}
		\end{equation}
		and
		\begin{equation}\label{Gamma-}
			\{\text{poles of $\psi_S(x)$ but not of $\psi_S(x) + \gamma  \psi_{S,T}(x)$}\}
		\end{equation}
	\end{enumerate}
	Moreover, if $S$ is $\gamma$-strongly-cospectral to $T$, then $\Lambda_{S,T}^{\gamma}$ is given by Equation \eqref{Gamma+}, and $\Lambda_{S,T}^{-\gamma}$ is given by Equation \eqref{Gamma-}.
\end{theorem}
\proof
Let $H$ be the Hermitian adjacency matrix of $Y$ with spectral decomposition
\[H = \sum_{\lambda} \lambda E_{\lambda}.\]
Suppose first that $S$ is $\gamma$-strongly-cospectral to $T$. Then $\abs{\gamma}=1$, and for each $\lambda$,
\begin{equation}\label{sc}
	E_{\lambda} B_S = \pm \gamma E_{\lambda} B_T.
\end{equation}
Taking matrix inner product on both sides of Equation \eqref{sc} gives
\[\tr((E_{\lambda})_{S,S})=\ip{E_{\lambda}B_S}{E_{\lambda}B_S} = \ip{E_{\lambda}B_T}{E_{\lambda}B_T} = \tr((E_{\lambda})_{T,T})\]
for each $\lambda$. This shows (i).  Now multiply both sides of Equation \eqref{sc} on the left by $B_S$ and take the trace. We obtain
\begin{equation}\label{tr_rln}
	\tr((E_{\lambda})_{S,S}) = \tr(B_S E_{\lambda} B_S )= \pm \gamma \tr(B_SE_{\lambda} B_T) = \pm \gamma \tr((E_{\lambda})_{S,T}) .
	\end{equation}
Thus
\begin{align*}
	\gamma \psi_{S,T}(x)&= \sum_{\lambda} \frac{\gamma }{x-\lambda} \tr( (E_{\lambda})_{S,T} )\\
	&=\sum_{\lambda\in \Lambda_{S,T}^+} \frac{1}{x-\lambda} \tr( (E_{\lambda})_{S,S} ) - \sum_{\lambda\in \Lambda_{S,T}^-} \frac{1}{x-\lambda}  \tr( (E_{\lambda})_{S,S} ),
\end{align*}
which, since $H$ is Hermitian, is a rational function over $\re$. This shows (ii). Finally, let $\lambda$ be a pole of $\psi_S(x)$. By Lemma \ref{evsp}, $\lambda$ lies in the eigenvalue support of $S$ and so $\tr((E_{\lambda})_{S,S})\ne 0$. Since
\[\psi_S(x) \pm \gamma \psi_{S,T}(x) = \sum_{\lambda} \frac{1}{x-\lambda}(\tr((E_{\lambda})_{S,S} \pm \gamma \tr((E_{\lambda})_{S,T} )),\]
by Equation \eqref{tr_rln}, either $\lambda$ is not a pole of $\psi_S(x)-\gamma \psi_{S,T}(x)$, or it is not a pole of $\psi_S(x)+\gamma \psi_{S,T}(x)$. This shows (iii).

Now suppose (i), (ii) and (iii) hold. By (i), for each $\lambda$ we have
\[
	\tr((E_{\lambda})_{S,S})  =  \tr((E_{\lambda})_{T,T}) .
\]
By (ii), there exists unimodular $\gamma$ such that
\[\gamma \psi_{S,T}(x) = \frac{p(x)}{q(x)}\]
for some coprime $p(x), q(x)\in \re[x]$. Since 
\begin{align*}
	\gamma \psi_{S,T}(x) &= \gamma \tr((xI-I)^{-1}_{S,T})=\sum_{\lambda} \frac{1}{x-\lambda} \gamma \tr((E_{\lambda})_{S,T}),
\end{align*}
we have
\[\gamma \tr((E_{\lambda})_{S,T}) = \lim_{x\to \lambda} (x-\lambda) \frac{p(x)}{q(x)}\]
which must be real. By (iii), for any $\lambda$ in the eigenvalue support of $S$, we have
\[
\tr((E_{\lambda})_{S,S})  = \pm \gamma \tr((E_{\lambda})_{S,T}) .
\]
Since 
\[\abs{\ip{E_{\lambda}B_S}{E_{\lambda}B_T}}^2 \le \ip{E_{\lambda}B_S}{E_{\lambda}B_S} \ip{E_{\lambda}B_T}{E_{\lambda}B_T},\]
with equality if and only if $E_{\lambda}B_S$ is a scalar multiple of $E_{\lambda} B_T$, it follows that 
\[E_{\lambda}B_S = \pm \gamma E_{\lambda} B_T.\]
Therefore, $S$ is $\gamma$-strongly-cospectral to $T$.
\qed

We can say more if $H$ is diagonally similar to a matrix with rational entries.

\begin{lemma}\label{rat_S_T}
Let $Y$ be a Hermitian weighted digraph with Hermitian adjacency matrix $H$. Let $S$ and $T$ be two subsets of vertices in $Y$ with equal size. Suppose there is a diagonal matrix $\Delta$ such that both $\Delta H \Delta^{-1}$ and $\Delta_{S,S} \Delta^{-1}_{T,T}$ have rational entries. Then $\psi_{S,T}(x)$ is a rational function over $\rats$.
\end{lemma}
\proof
We can write 
\[\psi_{S,T}(x) = \frac{\tr(\adj(xI-H)_{S,T})}{\det(xI-H)} = \frac{p(x)}{q(x)}\]
for some $p(x), q(x)\in \cx[x]$. Since $H$ is diagonally similar to a matrix with rational entries, the denominator $q(x)$ lies in $\rats[x]$. Thus it suffices to show that $p(x)\in \rats[x]$.

Let $M= \Delta H \Delta^{-1}$. Order the elements of $S$ and $T$ as they appear in the indices of $H$, say $S=\{\seq{a}{1}{2}{d}\}$ and $T=\{\seq{b}{1}{2}{d}\}$. For $j=1,2,\cdots,d$, we have
\begin{align*}
	\adj(xI-H)_{a_j, b_j}&= (-1)^{a_j+b_j} \det((xI-H)\backslash[b_j, a_j])\\
	&=(-1)^{a_j+b_j}\det((\Delta^{-1}(xI-M) \Delta)\backslash [b_j, a_j])\\
	&=(-1)^{a_j+b_j}\det(\Delta^{-1}\backslash [b_j, b_j])\det( (xI-M)\backslash [b_j, a_j] )\det(\Delta\backslash [a_j, a_j]).
\end{align*}
Since $\Delta_{S,S}\Delta^{-1}_{T,T}$ has rational entries, 
\[\det(\Delta^{-1}\backslash [b_j, b_j])\det(\Delta\backslash [a_j, a_j]) = \Delta_{a_j,a_j}\Delta^{-1}_{b_j, b_j}\in \rats,\]
and so we have 
\[\adj(xI-H)_{a_j, b_j} \in \rats[x].\]
Therefore
\[	p(x) = \sum_{j=1}^d \adj(xI-H)_{a_j, b_j} \in \rats[x].\tag*{\sqr53}\]

Combining this with Theorem \ref{gamma-sc-char} yields the following observation on the eigenvalue supports of $S$ and $T$.

\begin{corollary}\label{rat_sc}
	Let $Y$ be a Hermitian weighted digraph with Hermitian adjacency matrix $H$. Let $S$ and $T$ be two subsets of vertices in $Y$ with equal size. Suppose there is a diagonal matrix $\Delta$ such that both $\Delta H \Delta^{-1}$ and $\Delta_{S,S} \Delta^{-1}_{T,T}$ have rational entries. If $S$ is $\gamma$-strongly-cospectral to $T$, then they are strongly-cospectral, and $\Lambda_{S,T}^+$ and $\Lambda_{S,T}^-$ are each closed under algebraic conjugates.
\end{corollary}
\proof
The fact that $\gamma=\pm 1$ follows from Theorem \ref{gamma-sc-char} (ii) and Lemma \ref{rat_S_T}. From the assumption on $H$, $S$ and $T$ we see that $\psi_{S}(x)$ is a rational function over $\rats$, and so is
\[\psi_S(x) \mp \psi_{S,T}(x).\]
Therefore the poles of $\psi_S(x) \mp \psi_{S,T}(x)$ are closed under algebraic conjugates. By Theorem \ref{gamma-sc-char}, $\Lambda_{S,T}^{\pm}$ is closed under algebraic conjugates.
\qed

\section{Rational subspaces}\label{ptwise_psst}
We say a subspace $W$ is \textsl{rational} if the orthogonal projection onto $W$ has rational entries, or equivalently, if it has a basis with rational entries. Given two rational subspaces $W\sbs \col(C_a+I)$ and $V\sbs \col(C_b+I)$, and a rational coin operator $C$, there is an orthogonal coin basis induced by $(W,V)$ with rational entries. Let $M$ be the matrix with this basis as columns and 
\[N = M (M^TM)^{-1/2}.\]
Then $C=2NN^T-I$. Moreover, 
\[N^TRN = (M^TM)^{-1/2}M^TRM (M^TM)^{-1/2},\]
which is diagonally similar to a matrix with rational entries. If $S$ and $T$ are the $W$-clones of $a$ and the $V$-clones of $b$, then
\[(M^TM)_{S,S} = (M^TM)_{T,T},\]
and so by Lemma \ref{rat_S_T}, $\psi_{S,T}(x)$ is a rational function over $\rats$.

The first consequence of this observation is that, given rational $W$ and $C$, pointwise $W$-periodicity with integer periods can be tested in polynomial time. This extends Theorem 3.2 in \cite{Heller2025} on one-dimensional subspace periodicity with Grover coins. To start, we introduce some tools from number theory. The \textsl{$m$-th cyclotomic polynomial} is given by 
\[\Phi_m(x)=\prod_{1\le j\le m,\;\gcd(j,m)=1} \left(x - e^{2\pi ij/m}\right).\]
Let $\varphi(\cdot)$ denote the Euler totient function. Then $\deg(\Phi_m(x)) = \varphi(m)$. The following result from \cite{Bradford2005}  gives a lower bound on $\varphi(m)$ in terms of $m$.

\begin{theorem}\cite{Bradford2005}\label{totient_bound}
	For all natural numbers $m\ge 2$, we have
	\[m\le 3\varphi(m)^{3/2}.\]
\end{theorem}

By Theorem \ref{char_psp}, pointwise $W$-periodicity with integer periods translates into a number theoretic condition on the eigenvalue support of the $W$-clones, which, due to Theorem \ref{totient_bound}, can be checked in polynomial time. For ease of notation, given any polynomial $h$,  let $h^{\sharp}$ denote the polynomial obtained by
\[h^{\sharp}(x) = 2^{\deg(h)} x^{\deg(h)} h\left( \frac{x+x^{-1}}{2}\right).\]

\begin{theorem}\label{poly-time-sp}
	Let $C$ be a rational coin operator and $W$ a rational subspace of $\col(C_a+I)$. The following algorithm tests whether $U$ is pointwise $W$-periodic at $a$ with integer periods and, if this happens, returns the minimum period. Moreover, the running time of this algorithm is polynomial in the number of vertices of the original graph.
	\begin{enumerate}[(1)]
		\item 	Construct a digraph $Y$ relative to any orthonormal coin basis induced by $W$, and let $S$ be the $W$-clones of $a$. 
		\item Write
		\[\psi_S(x) = \frac{p(x)}{q(x)}\]
		for some coprime $p(x), q(x) \in \rats[x]$.
		\item  Let 
		\[g(x) = \frac{q(x)}{\gcd(p(x), q(x))}.\]
		\item 
			If
		\begin{equation*}\label{cyc_factor}
			g^{\sharp}(x)= \prod_{m\in L} \Phi_{m}(x),
		\end{equation*}
		for some set $L$ of positive integers, then $U$ is pointwise $W$-periodic at $a$ with minimum period $\lcm\{m: m\in L\}$. Otherwise, $U$ is not pointwise $W$-periodic at $a$ with integer periods.
	\end{enumerate}
\end{theorem}
\proof
By Theorem \ref{char_psp}, $U$ is pointwise $W$-periodic with integer period $\tau$ if and only if for each $\lambda\in \Lambda_S$,
\[\lambda = \frac{1}{2}(e^{j\pi/\tau} + e^{-j\pi/\tau})\]
for some even integer $j$, or by Lemma \ref{evsp}, the polynomial $h^{\sharp}(x)$ divides $x^{\tau}-1$. Thus $U$ is pointwise $W$-periodic at $a$ if and only if $h^{\sharp}(x)$ factors into cyclotomic polynomials, and the minimum period is given by the least common multiplier of the orders of these cyclotomic polynomials.

Let $X$ be the original graph with $n$ vertices. Since each vertex $u$ in $X$ has at most $\deg(u)$ clones, the digraph $Y$ constructed in Step (1) has $O(n^2)$ vertices. Since $\psi_S(x)$ is a rational function over $\rats$, computing $p(x),q(x)$ in Step (2) and $g(x)$ in Step (3) takes time polynomial in $n$. Since $g^{\sharp}(x)$ lies in $\rats[x]$, factoring it into irreducible divisors in $\rats[x]$ takes time polynomial in $n$, and by Theorem \ref{totient_bound}, it suffices to compare them with cyclotomic polynomials $\Phi_m(x)$ with $m$ up to $3n^3$ in Step (4). 
\qed

We are primarily interested in pointwise perfect $W$-transfer at integer steps where both $C$ and $W$ are rational. Using Theorem \ref{poly-time-sp} and results from Section \ref{gamma-sc},, we show this can be tested in polynomial time.
		
		\begin{theorem}
			Let $C$ be a rational coin and $W$ a rational subspace of $\col(C_a+I)\cap \col(C_b+I)$. The following algorithm tests whether pointwise perfect $W$-transfer from $a$ to $b$ occurs at integer steps and, if this happens, returns the minimum step needed for such transfer. Moreover, the running time of this algorithm is polynomial in the number of vertices of the original graph.
			\begin{enumerate}[(1)]
				\item Construct a digraph $Y$ relative to any orthonormal coin basis induced by $W$, and let $S$ and $T$ be the $W$-clones of $a$ and the $W$-clones of $b$, respectively.
				\item Check if $\psi_S(x) = \psi_T(x)$. If not, pointwise perfect $W$-transfer does not occur from $a$ to $b$.
				\item Write 	\[\psi_S(x) = \frac{p(x)}{q(x)}\]
				for some coprime $p(x), q(x) \in \rats[x]$. Let 
				\[g(x) = \frac{q(x)}{\gcd(p(x), q(x))}.\]
				\item Check if $g^{\sharp}(x)$ factors into cyclotomic polynomials. If not, pointwise perfect $W$-transfer does not occur from $a$ to $b$ at any integer step. Otherwise, let $L$ be the set of positive integers such that
				\[
					g^{\sharp}(x)= \prod_{m\in L} \Phi_{m}(x),
			\]
			and let $\tau =\lcm\{m: m \in L\}$. 
			\item Check if $\tau$ is even. If not, pointwise perfect $W$-transfer does not occur from $a$ to $b$ at any integer step. Otherwise, let
			\[L_+ = \{m\in L: \tau /m \text{ is even}\},\quad L_- = \{m\in L: \tau/m \text{ is odd}\}.\]
			\item Write 
			\[\psi_S(x) - \psi_{S,T}(x) = \frac{p_+(x)}{q_+(x)},\quad \psi_S(x) + \psi_{S,T}(x) = \frac{p_-(x)}{q_-(x)}\]
			for some  $p_+(x), q_+(x), p_-(x), q_-(x) \in \rats[x]$. Let 
			\[g_+(x) = \frac{q_+(x)}{\gcd(p_+(x), q_+(x))},\quad g_-(x) = \frac{q_-(x)}{\gcd(p_-(x), q_-(x))}.\]
			Check if 
			\[g^{\sharp}_+(x)= \prod_{m\in L_+} \Phi_{m}(x),\quad g^{\sharp}_-(x) = \prod_{m\in L_-} \Phi_{m}(x),\]
			or
			\[g^{\sharp}_+(x)= \prod_{m\in L_-} \Phi_{m}(x),\quad g^{\sharp}_-(x) = \prod_{m\in L_+} \Phi_{m}(x).\]
			If so, pointwise perfect $W$-transfer occurs from $a$ to $b$ at time $\tau/2$. Otherwise, it does not occur at any integer step.
			\end{enumerate}
			\end{theorem}
			\proof
			By Theorem \ref{poly-time-sp}, the quantum walk passes Step (3) and (4) if and only if $U$ is $W$-periodic at $a$ with mininum integer period $\tau$. Thus the minimum time for pointwise perfect $W$-transfer to occur from $a$ to $b$ is $\tau/2$, which is an integer if and only if $\tau$ is even.
			
			Now suppose $U$ is pointwise $W$-periodic at $a$ with even minimum period $\tau$. Let $t=\tau/2$. By Theorem \ref{char-ptwise-psst} and Corollary \ref{rat_S_T}, it suffices to show that Steps (2), (5) and (6) test strong cospectrality between $S$ and $T$ and the parity condition for $\Lambda_{S,T}^{\pm}$ in polynomial time. Indeed, as $C$ and $W$ are both rational, $S$ and $T$ are strongly cospectral if and only if they are cospectral and each irreducible factor of $g^{\sharp}(x)$ divides exactly one of  $g^{\sharp}_+(x)$ and $g^{\sharp}_-(x)$. Note that if $\tau/m$ is even, then all primitive $m$-th roots unity are also $t$-th roots of unity, and if $\tau/m$ is odd, then no primitive $m$-th root of unity is a $t$-th root of unity. Therefore, we have pointwise perfect $W$-transfer in the form
			\[U^t x_a(w) = x_b(w),\quad \forall w\in W\]
			if and only if 
			\[g^{\sharp}_+(x)= \prod_{m\in L_+} \Phi_{m}(x),\quad g^{\sharp}_-(x) = \prod_{m\in L_-} \Phi_{m}(x),\]
			and pointwise perfect $W$-transfer in the form
			\[U^t x_a(w) = -x_b(w),\quad \forall w\in W\]
			if and only if 
			\[g^{\sharp}_+(x)= \prod_{m\in L_-} \Phi_{m}(x),\quad g^{\sharp}_-(x) = \prod_{m\in L_+} \Phi_{m}(x). \tag*{\sqr53}\]
			
\section{Transfer between marked vertices}			
Working with $H$ instead of $U$ reduces the scale of the problem, especially when most coins are reflections about low dimensional subspaces. In this section, we consider a special arc-reversal walk where $C$ assigns the Grover coin
\[G_u = \frac{2}{\deg(u)}J-I\]
to any vertex $u$ other than $a$ and $b$. As the Grover coin is a reflection about $\col(\one)$, our reduction yields a weighted Hermitian digraph with size close to the original graph $X$. The exact structure of the digraph depends on the reflection coins $C_a$ and $C_b$ we assign to the $a$ and $b$. Following the notion of quantum walk search algorithms (see, for example, \cite{Shenvi2003,Wong2018,Stefanak2016}), we will call such a walk a Grover coined walk with \textsl{marked vertices} $a$ and $b$.

Write 
\[C_a = 2KK^*-I,\quad C_b = 2LL^*-I\]
for some matrices $K$ and $L$ with orthonormal columns. Let $Y$ be the digraph relative to an orthonormal coin basis containing columns of $K$ and columns of $L$. Depending on whether $a$ and $b$ are adjacent in $X$, there are two possible underlying graphs of $Y$, as illustrated in Figure \ref{Y}. Here $k$ and $l$ are the number of columns of $K$ and $L$, respectively.

\begin{figure}[h]
	\centering
	\begin{minipage}[t]{0.5\textwidth}
		\centering
		\begin{tikzpicture}
			\draw(0,0) ellipse [x radius=2.5, y radius=1.8];
			\draw [color=blue!60,very thick][rotate=28] (-1.7,0.3,-2.2)
			ellipse [x radius=1.3, y radius=0.7];
			\node [rotate=40][scale=1, color=blue] at (-1.6,0.4) {$N(a)$};
			\draw [thick,color=blue!60, very thick](-2.5,0.2) -- (-2.9,1.2);
			\draw [thick,color=blue!60, very thick](-0.15,1.25) -- (-2.9,1.2);
			\draw [thick,color=blue!60, very thick](-2.5,0.2) -- (-1.5,2);
			\draw [thick,color=blue!60, very thick](-0.15,1.25) -- (-1.5,2);
			\node [rotate=30][scale=1.2, color=blue] at (-2.3,2.0) {$a_1,\dots,a_{k}$};

			\draw [color=red!60,very thick][rotate=-5](-0.6,0.25,-2.2)
			ellipse [x radius=1.3, y radius=0.7];
			\node [scale=1, color=red] at (0.6,1.0) {$N(b)$};
			\draw [thick,color=red!60, very thick](-0.93,1.25) -- (-0.45,2.3);
			\draw [thick,color=red!60, very thick](1.6,1.15) -- (-0.45,2.3);
			\draw [thick,color=red!60, very thick](-0.93,1.25) -- (1.35,2.3);
			\draw [thick,color=red!60, very thick](1.6,1.15) -- (1.35,2.3);
			\node [scale=1.2, color=red] at (0.6,2.6) {$b_1,\dots,b_{\ell}$};
			
			\node [scale=1.1] at (0,-0.6) {$X\setminus\{a,b\}$};
		\end{tikzpicture}
		\caption*{If  $a\not\sim b$ in $X$}
	\end{minipage}%
	\begin{minipage}[t]{0.5\textwidth}
		\centering
		\begin{tikzpicture}
			\draw(0,0) ellipse [x radius=2.5, y radius=1.8];
			\draw [color=blue!60,very thick][rotate=28] (-1.7,0.3,-2.2)
			ellipse [x radius=1.3, y radius=0.7];
			\node [rotate=33][scale=1, color=blue,inner sep=1pt] at (-1.6,0.4) {$N(a)\setminus b$};
			\draw [color=blue!60, very thick](-2.5,0.2) -- (-2.9,1.2);
			\draw [color=blue!60, very thick](-0.15,1.25) -- (-2.9,1.2);
			\draw [color=blue!60, very thick](-2.5,0.2) -- (-1.5,2);
			\draw [color=blue!60, very thick](-0.15,1.25) -- (-1.5,2);
			\node [rotate=30][scale=1.1, color=blue] at (-2.3,2.0) {$a_1,\dots,a_k$};

			\draw [color=red!60,very thick][rotate=-5](-0.6,0.25,-2.2)
			ellipse [x radius=1.3, y radius=0.7];
			\node [scale=1, color=red,inner sep=1pt] at (0.75,1) {$N(b)\setminus a$};
			\draw [color=red!60, very thick](-0.93,1.25) -- (-0.45,2.3);
			\draw [color=red!60, very thick](1.6,1.15) -- (-0.45,2.3);
			\draw [color=red!60, very thick](-0.93,1.25) -- (1.35,2.3);
			\draw [color=red!60, very thick](1.6,1.15) -- (1.35,2.3);
			\node [scale=1.1, color=red] at (0.6,2.6) {$b_1,\dots,b_{\ell}$};
			\node [scale=1.1] at (0,-0.6) {$X\setminus\{a,b\}$};
			
			\draw(-1.2,2.5)[color=blue!60,very thick, rotate=25] arc [start angle=-355,
			end angle=-185, x radius=1.25, y radius=0.6];
			\draw(1.8,2.62)[color=red!60,very thick] arc [start angle=-355,
			end angle=-185, x radius=1.25, y radius=0.6];
			\draw(0.5,3.2)[rotate=13] arc [start angle=-355,
			end angle=-185, x radius=1.6, y radius=0.7];
		\end{tikzpicture}
		\caption*{If $a\sim b$ in $X$}
	\end{minipage}
	\caption{The underlying graph of $Y$}
	\label{Y}
\end{figure}
		
If $a$ and $b$ are adjacent in $X$, then perfect subspace state transfer from $a$ to $b$ is guaranteed at $t=1$:
\[U (Ce_{(a,b)}) = R e_{(a,b)} = e_{(b,a)}.\]
Hence from this point on, we will only consider subspace transfer between non-adjacent vertices. Let $\Delta$ be the diagonal matrix obtained from the degree matrix $D$ of $X$ by removing the rows and columns indexed by $a$ and $b$:
\[\Delta = D\setminus \{a,b\}.\] 
Given $a$ and $b$ are non-adjacent, we can write the Hermitian adjacency matrix $H$ of $Y$ as follows:
\begin{center}
\begin{tikzpicture}
	\draw [thick](0.25,0) -- (0.25,5.5);
	\draw [thick](0.25,0) -- (0.5,0);
	\draw [thick](0.25,5.5) -- (0.5,5.5);
	\draw [thick](7,0) -- (7,5.5);
	\draw [thick](6.75,0) -- (7,0);
	\draw [thick](6.75,5.5) -- (7,5.5);
	\draw [thick] (3,0) rectangle (6.5,2.75);
	\node [scale=1] at (4.75,1.375) {$A(X\setminus\{a,b\})$};

	\draw [thick,color=blue!60] (0.5,1.5) rectangle (1.5,3);
	\node [scale=1,color=blue] at (1,2.25) {$K$};
	\draw [thick,color=red!60] (1.75,0.75) rectangle (2.75,2.25);
	\node [scale=1,color=red] at (2.25,1.5) {$L$};
	
	\draw [thick,color=blue!60] (3,4.25) rectangle (4.5,5.25);
	\node [scale=1,color=blue] at (3.75,4.725) {$K^*$};
	\draw [thick,color=red!60] (4,3) rectangle (5.5,4);
	\node [scale=1,color=red] at (4.75,3.5) {$L^*$};

	\node [scale=1] at (1,4.75) {$0$};
	\node [scale=1] at (2.25,4.75) {$0$};
	\node [scale=1] at (1,3.5) {$0$};
	\node [scale=1] at (2.25,3.5) {$0$};
	
	\node [scale=0.7,color=blue] at (1,5.75) {$a_1,\cdots, a_k$};
	\node [scale=0.7,color=red] at (2.25,5.75) {$b_1,\cdots,b_{\ell}$};
	
	\node [scale=0.7,color=blue] at (3.75,5.75) {$N(a)$};
	\draw [color=blue!60](3,5.5) -- (4.5,5.5);
	\node [scale=0.7,color=red] at (4.75,5.75) {$N(b)$};
	\draw [color=red!60](4,5.45) -- (5.5,5.45);
	
	\node [scale=1] at (-1.25,5.25) {$H=\pmat{I & \\ & \Delta}^{-\frac{1}{2}}$};
	\node [scale=1] at (8,5.25) {$\pmat{I & \\ & \Delta}^{-\frac{1}{2}}$};
\end{tikzpicture}
\end{center}

Since most coins in this walk are reflections about one-dimensional subspaces, we can afford a characterization of pointwise perfect subspace state transfer between $a$ and $b$ that involves only  the coins $C_a$ and $C_b$ rather than their specific decompositions into $K$ and $L$. To start, let  $M$ be the block diagonal matrix
\[M = \pmat{K & & \\ & L & \\ & & I}\]
where the last identity matrix has the same size as  $\Delta$. 
Then $MM^*=I$ and 
\[MHM^* = G\]
where
\begin{center}
\begin{tikzpicture}
	\draw [thick](0.25,0) -- (0.25,5.5);
	\draw [thick](0.25,0) -- (0.5,0);
	\draw [thick](0.25,5.5) -- (0.5,5.5);
	\draw [thick](7,0) -- (7,5.5);
	\draw [thick](6.75,0) -- (7,0);
	\draw [thick](6.75,5.5) -- (7,5.5);
	\draw [thick] (3,0) rectangle (6.5,2.75);
	\node [scale=1] at (4.75,1.375) {$A(X\setminus\{a,b\})$};

	\draw [thick,color=blue!60] (0.5,1.5) rectangle (1.5,2.5);
	\node [scale=1,color=blue] at (1,2) {$\frac{C_a+I}{2}$};
	\draw [thick,color=red!60] (1.75,1) rectangle (2.75,2);
	\node [scale=1,color=red] at (2.25,1.5) {$\frac{C_b+I}{2}$};
	
	\draw [thick,color=blue!60] (3,4.25) rectangle (4,5.25);
	\node [scale=1,color=blue] at (3.5,4.725) {$\frac{C_a+I}{2}$};
	\draw [thick,color=red!60] (3.5,3) rectangle (4.5,4);
	\node [scale=1,color=red] at (4,3.5) {$\frac{C_b+I}{2}$};

	\node [scale=1] at (2.25,0.5) {$0$};
	\node [scale=1] at (1,0.5) {$0$};
	
	\node [scale=1] at (5.5,4.725) {$0$};
	\node [scale=1] at (5.5,3.5) {$0$};
	
	\node [scale=1] at (1,4.75) {$0$};
	\node [scale=1] at (2.25,4.75) {$0$};
	\node [scale=1] at (1,3.5) {$0$};
	\node [scale=1] at (2.25,3.5) {$0$};
	
	\node [scale=0.7,color=blue] at (1,5.75) {$\cl(a)$};
	\node [scale=0.7,color=red] at (2.25,5.75) {$\cl(b)$};
	
	\node [scale=0.7,color=blue] at (3.4,5.75) {$N(a)$};
	\draw [color=blue!60](3,5.5) -- (4,5.5);
	\node [scale=0.7,color=red] at (4.1,5.75) {$N(b)$};
	\draw [color=red!60](3.5,5.45) -- (4.5,5.45);
	
	\node [scale=1] at (-1.2,5.25) {$G=\pmat{I & \\ & \Delta}^{-\frac{1}{2}}$};
	\node [scale=1] at (8,5.25) {$\pmat{I & \\ & \Delta}^{-\frac{1}{2}}$};
\end{tikzpicture}
\end{center}
The underlying digraph of $G$ can be seen as obtained from $X$ by ``blowing up" $a$ and $b$ into $\deg(a)$ clones of $a$, denoted $\cl(a)$,  and $\deg(b)$ clones of $b$, denoted $\cl(b)$, whose connections to the rest of the graph is described by $C_a$ and $C_b$. We will refer to such digraph as the \textsl{$(C_a,C_b)$-blow-up of $X$}. Here we kept the name \textsl{clones}  with the understanding that they arise from the projections onto $\col(C_a+I)$ and $\col(C_b+I)$ rather than specific bases for these subspaces.

Note that $\col(M)$ is $G$-invariant. Thus if
\[G=\sum_{\lambda} \lambda E_{\lambda}\]
is the spectral decomposition of $G$, then
\[H = \sum_{\lambda: E_{\lambda}M\ne 0} \lambda (M^* E_{\lambda} M)\]
is the spectral decomposition of $H$. As a result, for any subspace $W$ of $\col(C_a+I)$ and the set $S$ of $W$-clones relative to $H$, the eigenvalue $\lambda$ of $H$ lies in $\Lambda_S$ if and only if the column space of $E_{\lambda}[\cl(a), \cl(a)]$ is not orthogonal to $W$. Since $G$ is independent of the basis for $W$, we will also refer to $\Lambda_S$ as the \textsl{eigenvalue support of $x_a(W)$}, and denote it by $\Lambda_{a(W)}$:
\[\Lambda_{a(W)} =  \{\lambda: 	 E_{\lambda}[\cl(a),\cl(a)]\upharpoonright_W\ne 0\}.\]
Moreover, given $W\sbs \col(C_a+I) \cap \col(C_b+I)$, the $W$-clones of $a$ are $\gamma$-strongly-cospectral to the $W$-clones of $b$ if and only if 
\[E_{\lambda}[\cl(a),\cl(a)\cup \cl(b)]\upharpoonright_W =\pm \gamma  E_{\lambda}[\cl(b),\cl(a)\cup \cl(b)]\upharpoonright_W\]
for each eigenvalue $\lambda$ of $G$. When this happens, we  say \textsl{$x_a(W)$ is $\gamma$-strongly cospectral to $x_b(V)$}, and define
\[\Lambda_{a(W), b(W)}^{\gamma}=\{\lambda: E_{\lambda}[\cl(a),\cl(a)\cup \cl(b)]\upharpoonright_W = \gamma  E_{\lambda}[\cl(b),\cl(a)\cup \cl(b)]\upharpoonright_W \}\]
and 
\[\Lambda_{a(W), b(W)}^{-\gamma}=\{\lambda: E_{\lambda}[\cl(a),\cl(a)\cup \cl(b)]\upharpoonright_W =- \gamma  E_{\lambda}[\cl(b),\cl(a)\cup \cl(b)]\upharpoonright_W \}.\]
This, together with Theorem \ref{char-ptwise-psst} and Theorem \ref{char-ptwise-pgsst}, leads to  basis-free characterizatiosn of pointwise perfect and pretty good $W$-transfer between marked vertices.

\begin{theorem}\label{basis-free}
	Let $U=RC$ be the transition matrix of a Grover coined walk on $X$ with marked vertices $a$ and $b$. Let $W\sbs \col(C_a + I)\cap \col(C_b+I)$ be a subspace. Let 
	\[G = \sum_{\lambda} \lambda E_{\lambda}\]
	be the spectral decomposition of the $(C_a, C_b)$-blow-up of $X$. There is pointwise perfect $W$-transfer from $a$ to $b$ at time $t$ if and only if the following hold for some unimodular $\gamma$.
	\begin{enumerate}[(i)]
		\item  $x_a(W)$ is $\gamma$-strongly-cospectral to $x_b(W)$.
		\item  $t\arccos\lambda$ is an even multiple of $\pi$ for each eigenvalue $\lambda\in \Lambda_{a(W),b(W)}^{\gamma}$.
		\item $t\arccos\lambda$ is an odd multiple of $\pi$ for each eigenvalue $\lambda\in \Lambda_{a(W),b(W)}^{-\gamma}$.
	\end{enumerate}
\end{theorem}		

\begin{theorem}\label{basis-free-pgst}
	Let $U=RC$ be the transition matrix of a Grover coined walk on $X$ with marked vertices $a$ and $b$. Let $W\sbs \col(C_a + I)\cap \col(C_b+I)$ be a subspace. Let 
	\[G = \sum_{\lambda} \lambda E_{\lambda}\]
	be the spectral decomposition of the $(C_a, C_b)$-blow-up of $X$. There is pointwise pretty good $W$-transfer from $a$ to $b$ if and only if the following hold for some unimodular $\gamma$.
	\begin{enumerate}[(i)]
		\item  $x_a(W)$ is $\gamma$-strongly-cospectral to $x_b(W)$.
				\item  For any set $\{\ell_{\lambda}:\lambda\in\Lambda_{a(W)}\}$ of integers such that
		\begin{equation}\label{K3}
			\sum_{\lambda \in \Lambda_{a(W)}} \ell_{\lambda} \arccos \lambda =0,
		\end{equation}
		we have
		\[\sum_{\lambda\in\Lambda^{-\gamma}_{a(W),b(W)}}\ell_{\lambda} \equiv 0\pmod2.\] 
		\end{enumerate}  
				Moreover, the target states can be approximated in integer steps if and only if (i) and (ii) hold with Equation \eqref{K3} replaced by
		\[\sum_{\lambda \in \Lambda_{a(W)}} \ell_{\lambda} \arccos \lambda \equiv 0 \pmod{2\pi}.\]
\end{theorem}		

When an explicit eigenbasis for $G$ is known, we may work with the following equivalent definitions of eigenvalue support and $\gamma$-strong-cospectrality.

\begin{lemma}
Let $U=RC$ be the transition matrix of a Grover coined walk on $X$ with marked vertices $a$ and $b$. Let $W\sbs \col(C_a + I)\cap \col(C_b+I)$ be a subspace.
	\begin{enumerate}[(i)]
		\item An eigenvalue $\lambda$ of $G$ lies in $\Lambda_{a(W)}$ if and only if some $\lambda$-eigenvector $z$ satisfies $\projonto{W} {z[\cl(a)]}\ne 0$.
		\item $x_a(W)$ is $\gamma$-strongly-cospectral to $x_b(W)$ if and only if for each eigenprojection $E_{\lambda}$ of $G$, there is $\sigma_{\lambda}\in \{1, -1\}$ such that every $\lambda$-eigenvector $z$ satisfies
		\[\projonto{W}{z[\cl(a)]} = \sigma_{\lambda} \gamma  \projonto{W}{z[\cl(b)]}.\]
	\end{enumerate}
\end{lemma}
\proof
The results follow from the fact that $E_{\lambda}$ is a sum of outer products of any orthonormal basis for the $\lambda$-eigenspace.
\qed

We conclude this section with an observation that, since $G$ is a block matrix, finding its eigenvectors  amounts to solving the quadratic eigenvalue problem in the following lemma with 

\begin{center}
\begin{tikzpicture}

	\draw[thick, color=blue!60] (0.5,1.5) rectangle (1.5,2.5);
	\node[scale=1, color=blue] at (1,2) {$\frac{C_a+I}{2}$};
	
	\draw[thick, color=red!60] (1.75,1) rectangle (2.75,2);
	\node[scale=1, color=red] at (2.25,1.5) {$\frac{C_b+I}{2}$};

	\node[scale=1] at (1,0.5) {$0$};
	\node[scale=1] at (2.25,0.5) {$0$};

	\draw[thick] (0.4,0) -- (0.4,2.6);
	\draw[thick] (0.4,0) -- (0.6,0);
	\draw[thick] (0.4,2.6) -- (0.6,2.6);

	\draw[thick] (2.85,0) -- (2.85,2.6);
	\draw[thick] (2.65,0) -- (2.85,0);
	\draw[thick] (2.65,2.6) -- (2.85,2.6);

	\node[scale=1] at (3.0,1.3) {,};
	
	% left label
	\node[scale=1] at (-0.6,1.55) {$F = \Delta^{-\frac{1}{2}}$};

	\node[scale=1, anchor=west] at (4.3,1.55)
	{$B = \Delta^{-\frac{1}{2}} A(X\setminus\{a,b\}) \Delta^{-\frac{1}{2}}$.};
\end{tikzpicture}
\end{center}

\begin{lemma}\label{quad_ev_eq}
	Let $G$ and $z$ be a block matrix and a block vector given by 
	\[G = \pmat{0 & F^* \\ F & B}, \quad z=\pmat{x\\ y}.\]
	Then  $z$ is an eigenvector for $G$ with eigenvalue $\lambda$ if and only if either 
	\[\lambda\ne0, \quad F^*y = \lambda x,\quad (\lambda^2I - \lambda B - FF^*) y=0\]
	or
	\[\lambda =0,\quad F^*y = 0,\quad Fx+By=0.\]
\end{lemma}
\proof
Expanding $Gz = \lambda z$ gives
\begin{equation}\label{qep}
	\begin{cases}
	F^*y = \lambda x,\\
	Fx + By = \lambda y.
\end{cases}
\end{equation}
If $\lambda=0$, the above is equivalent to
\[F^*y=0,\quad Fx + By =0.\]
Otherwise, multiplying the second equation by $\lambda$ and plugging the first equation in \eqref{qep} yields
\[ (\lambda^2I - \lambda B - FF^*) y=0.\]
Conversely, if 
\[\lambda\ne0, \quad F^*y = \lambda x,\quad (\lambda^2I - \lambda B - FF^*) y=0,\]
then it is not hard to see that both equations in \eqref{qep} hold.
\qed

\section{Infinite families}\label{inf_fam}

In this section, we construct infinite families of Grover coined walks that admit pointwise perfect or pretty good $d$-dimensional subspace state transfer between the marked vertices, where $d\ge 2$. We will assume $C_a = C_b$. Let $W$ be a subspace of $\col(C_a+I)$. We consider two special cases that make the analysis easier.

The first case is when $a$ and $b$ are twins, that is, when $N(a)=N(b)$. In this case, $a$ and $b$ are not adjacent, and we can write
\begin{center}
\begin{tikzpicture}
	\draw [thick](0.25,0) -- (0.25,5.5);
	\draw [thick](0.25,0) -- (0.5,0);
	\draw [thick](0.25,5.5) -- (0.5,5.5);
	\draw [thick](7,0) -- (7,5.5);
	\draw [thick](6.75,0) -- (7,0);
	\draw [thick](6.75,5.5) -- (7,5.5);
	\draw [thick] (3,0) rectangle (6.5,2.75);
	\node [scale=1] at (4.75,1.375) {$A(X\setminus\{a,b\})$};

	\draw [thick,color=blue!60] (0.5,1) rectangle (1.5,2);
	\node [scale=1,color=blue] at (1,1.5) {$\frac{C_a+I}{2}$};
	\draw [thick,color=red!60] (1.75,1) rectangle (2.75,2);
	\node [scale=1,color=red] at (2.25,1.5) {$\frac{C_b+I}{2}$};
	
	\draw [thick,color=blue!60] (3,4.25) rectangle (4.1,5.25);
	\node [scale=1,color=blue] at (3.55,4.725) {$\frac{C_a+I}{2}$};
	\draw [thick,color=red!60] (3,3) rectangle (4.1,4);
	\node [scale=1,color=red] at (3.55,3.5) {$\frac{C_b+I}{2}$};

	\node [scale=1] at (2.25,0.5) {$0$};
	\node [scale=1] at (1,0.5) {$0$};
	
	\node [scale=1] at (5.5,4.725) {$0$};
	\node [scale=1] at (5.5,3.5) {$0$};
	
	\node [scale=1] at (1,4.75) {$0$};
	\node [scale=1] at (2.25,4.75) {$0$};
	\node [scale=1] at (1,3.5) {$0$};
	\node [scale=1] at (2.25,3.5) {$0$};
	
	\node [scale=0.7,color=blue] at (1,5.75) {$\cl(a)$};
	\node [scale=0.7,color=red] at (2.25,5.75) {$\cl(b)$};
	
	\node [scale=0.7,color=blue] at (3.3,5.75) {$N(a)$};
	\draw [color=blue!60](3,5.5) -- (4.1,5.5);
	\node [scale=0.7,color=red] at (3.9,5.75) {$N(b)$};
	\draw [color=red!60](3,5.45) -- (4.1,5.45);
	
	\node [scale=1] at (-1.2,5.25) {$G=\pmat{I & \\ & \Delta}^{-\frac{1}{2}}$};
	\node [scale=1] at (8,5.25) {$\pmat{I & \\ & \Delta}^{-\frac{1}{2}}$};
\end{tikzpicture}

\end{center}
Note that if $C_a = C_b$, then $0$ is always an eigenvalue of $G$, and its eigenspace contains
\begin{equation}\label{0es_part}
	e_{1}-e_{\deg(a)+1}, \;e_{2}-e_{\deg(a)+2},\;\cdots, \;e_{\deg(a)} - e_{2\deg(a)}.
\end{equation}
Therefore, any eigenvector $z$ for $G$ orthogonal to these vectors must satisfy $z[\cl(a)] = z[\cl(b)]$. It follows that $x_a(W)$ is $\gamma$-strongly-cospectral to $x_b(W)$ if and only if it is strongly-cospectral to $x_b(W)$, which holds precisely when every vector  $z$ in $\ker(G)$ satisfies
\[\projonto{W}{z[\cl(a)]} = -\projonto{W}{z[\cl(b)]}.\]

For ease of notation, partition $A(X\setminus \{a,b\})$ as
\[A(X\setminus \{a,b\}) = \pmat{A_1 & A_2\\ A_2^T & A_3},\]
where $A_1$ is indexed by $N(a)$, and $A_3$ is indexed by the other vertices in $X\setminus\{a,b\}$.

\begin{theorem}\label{char-twins}
Let $U=RC$ be the transition matrix of a Grover coined walk on $X$ with marked vertices $a$ and $b$. Suppose $a$ and $b$ are twins and $C_a=C_b$. Let $W$ be a subspace of $\col(C_a + I)$.  Let 
	\[G = \sum_{\lambda} \lambda E_{\lambda}\]
	be the spectral decomposition of the $(C_a, C_b)$-expansion of $X$. Pointwise perfect $W$-transfer occurs from $a$ to $b$ at time $t$ if and only if $W$ is orthogonal to the column space of 
	\[\pmat{A_1 & A_2} \Delta^{-1/2} E_0[V(X)\setminus \{a,b\}, V(X)\setminus \{a,b\}],\]
	and one of the following holds:
	\begin{enumerate}[(i)]
		\item $t$ is a multiple of $4$, and $t \arccos\lambda$ 
		is an odd multiple of $\pi$ for each $\lambda\in \Lambda_{a(W)}\setminus\{0\}$.
		\item $t\equiv 2\pmod{4}$, and $t \arccos\lambda$ 
		is an even multiple of $\pi$ for each $\lambda\in \Lambda_{a(W)}\setminus\{0\}$.
	\end{enumerate}
\end{theorem}
\proof
Let $z$ be any $0$-eigenvector for $G$. Let $x$ and $y$ be the following components of $z$:
\[x = z[\cl(a)\cup  \cl(b)],\quad y = z[V(X)\setminus\{a,b\}].\]
By Lemma  \ref{quad_ev_eq}, 

\[
\pmat{
	\frac{C_a + I}{2} & \frac{C_b + I}{2} \\
	0 & 0
} x
= -A(X \setminus \{a,b\})\Delta^{-\frac{1}{2}} y.
\]
Let $P_W$ be the orthogonal projection onto $W$. Since $C_a=C_b$ and $W$ is a subspace of $\col(C_a+I)$,
\begin{align*}
 P_W{x[\cl(a)]} + P_W{x[\cl(b)]}  &= \pmat{P_W & 0} \pmat{
 	\frac{C_a + I}{2} & \frac{C_b + I}{2} \\
 	0 & 0
 } x\\
 &=-\pmat{P_W & 0}A(X \setminus \{a,b\})\Delta^{-\frac{1}{2}} y\\
&= -P_W{ \pmat{A_1 & A_2} \Delta^{-\frac{1}{2}} y}.
	\end{align*}
It follows from our earlier discussion that $x_a(W)$ and $x_b(W)$ are strongly-cospectral if and only if the $0$-eigenprojection $E_0$ of $G$ satisfies
\[P_W \pmat{A_1 & A_2} \Delta^{-\frac{1}{2}} E_0[V(X)\setminus\{a,b\}, V(X)\setminus\{a,b\}]=0.\]
In this case, $0$ is the only eigenvalue in $\Lambda_{a(W), b(W)}^-$, and since $\arccos 0 = \pi/2$, pointwise perfect $W$-transfer occurs from $a$ to $b$ if and only if one of the parity conditions in (i) and (ii) holds.
\qed

We define the \textsl{support} of a subspace $W$ to be
\[\supp(W) = \bigcup_{w\in W} \supp(w).\]
The following result describes $\Lambda_{a(W), b(W)}^+ $ and $\Lambda_{a(W), b(W)}^-$ for strongly-cospectral $x_a(W)$ and $x_b(W)$ when all vertices in $\supp(W)$ have the same degree.

\begin{corollary}\label{twin-supp}
	Let $U=RC$ be the transition matrix of a Grover coined walk on $X$ with marked vertices $a$ and $b$. Suppose $a$ and $b$ are twins and $C_a=C_b$. Let 
	\[W \sbs \col(C_a + I) \cap \ker\pmat{A_1 \\ A_2^T}.\]
	Then $x_a(W)$ and $x_b(W)$ are strongly-cospectral. Moreover, if $W$ is supported by vertices of the same degree $\delta$, then 
	\[ \Lambda_{a(W), b(W)}^+ =\left\{\pm \sqrt{\frac{2}{\delta}}\right\}, \quad \Lambda_{a(W), b(W)}^- = \{0\}.\]
\end{corollary}
\proof
Since $W$ is orthogonal to the column space of $\pmat{A_1 & A_2}$, by Theorem \ref{char-twins}, $x_a(W)$ and $x_b(W)$ are strongly-cospectral, and $\Lambda_{a(W), b(W)}^- = \{0\}$.

Let $P$ and $P_W$ be the projections onto $\col(C_a+I)$ and $W$, respectively. Let $\lambda$ be any nonzero eigenvalue of $G$. From Lemma \ref{quad_ev_eq} we see that 
\[G\pmat{x\\ y} = \lambda \pmat{x\\y}\]
if and only if
\begin{equation}\label{quad1}
	\pmat{P & 0\\ P&0} \Delta^{-1/2}y = \lambda x
\end{equation}
and
\begin{equation}\label{quad2}
\left(\lambda^2I - \lambda \Delta^{-1/2} A(X\setminus\{a,b\}) \Delta^{-1/2} - \Delta^{-1/2} \pmat{2P & 0\\0 &0} \Delta^{-1/2}\right) y=0.
\end{equation}
Pick any nonzero $w\in W$. Since all vertices in $\supp(W)$ have  degree $\delta$, 
\[\Delta^{-1/2} \pmat{w\\ 0} = \delta^{-1/2} \pmat{w\\ 0}.\]
Set $y=\pmat{w\\0}$ and plug it into the left-hand side of Equation \eqref{quad2}. We get
\begin{align*}
	&\left(\lambda^2I - \lambda \Delta^{-1/2} A(X\setminus\{a,b\}) \Delta^{-1/2} - \Delta^{-1/2} \pmat{2P & 0\\0 &0} \Delta^{-1/2}\right) y\\
	=&\lambda^2 \pmat{w\\ 0} -\lambda  \delta^{-1/2}\Delta^{-1/2} \pmat{A_1 \\ A_2^T } w - 2\delta^{-1}  \pmat{w\\ 0}\\
	=&(\lambda^2 - 2\delta^{-1})\pmat{w\\0},
\end{align*}
which is zero if and only if $\lambda=\pm \sqrt{2}/\sqrt{\delta}$. By Equation \eqref{quad1}, 
\[\pmat{\pm w/\sqrt{2}\\ \pm w/\sqrt{2} \\ w \\ 0}\]
is an eigenvector for $G$ with eigenvalue $\lambda=\pm \sqrt{2}/\sqrt{\delta}$. Thus for any $\mu$-eigenvector $z$ with $\mu\ne \pm \sqrt{2}/\sqrt{\delta}$, the component $z[N(a)]$ must be orthogonal to $W$. It follows from Equation \eqref{quad1} that $\mu$ does not lie in $\Lambda_{a(W), b(W)}^+$. Therefore, $\Lambda_{a(W), b(W)}^+=\{\pm \sqrt{2}/\sqrt{\delta}\}$.
\qed

We use this observation to construct three infinite families of walks that admit pointwise perfect subspace state transfer. The first family is able to transfer any subspace fixed by any reflection coin assigned to the marked vertices on $K_{2,m}$.

\begin{theorem}

Let $a$ and $b$ be the two vertices in $K_{2,m}$ of degree $m$. Let $U=RC$ be the transition matrix of a Grover coined walk on $K_{2,m}$ with marked vertices $a$ and $b$ and any reflection coins $C_a=C_b$. For any subspace $W$ of $\col(C_a+I)$, there is pointwise perfect $W$-transfer from $a$ to $b$ at time $t=2$.
\end{theorem}
\proof
Note that $A\setminus\{a,b\}=0$ and the zero subspace is orthogonal to any subspace $W$ of $\col(C_a+I)$. Since all neighbors of $a$ have degree $2$, by Corollary \ref{twin-supp}, $x_a(W)$ and $x_b(W)$ are strongly-cospectral with
\[\Lambda^+_{a(W), b(W)}=\{\pm1\} = \{\cos(0\pi/2),\cos(2\pi/2)\}\]
and 
\[ \Lambda^-_{a(W), b(W)}=\{0\} = \cos(\pi/2).\]
It follows from Theorem \ref{char-twins} (ii) that pointwise perfect $W$-transfer occurs from $a$ to $b$ at time $2$.
\qed

The next family is able to transfer a $2$-dimensional subspace between the marked vertices on certain $4$-regular circulant graphs.

\begin{theorem}\label{circulant}
For any positive integer $m$, let $c$ and $d$ be two distinct integers between $0$ and $m$ such that $c+d=m$. Let $U=RC$ be the transition matrix of a Grover coined walk on $X(2m, \pm \{c,d\})$ with marked vertices $0$ and $m$. 
Let $W$ be the subspace of $\cx^{\{c,d,-d,-c\}}$ spanned by
\[\left\{\pmat{1\\0\\-1\\0}, \pmat{0\\1\\0\\-1}\right\}.\]
Let $C_0=C_m$ be the reflection about any subspace containing $W$. There is pointwise perfect $W$-transfer from $0$ to $m$ at time $t=4$.
\end{theorem}
\proof
Let $a=0$ and $b=m$. First note that $a$ and $b$ are twins in $X$ as $c+d=m$:
\[N(a) =\{c, d, -d, -c\}=N(b).\]
Next, since $2c=-2d$ in $\ints_{2m}$, the vertices $c$ and $-d$ are twins in $X\setminus\{a,b\}$:
\[N(c) = \{2c, m, c-d, 0\} = N(-d)\]
and so are the vertices $-c$ and $d$:
\[N(d) = \{m, 2d, 0, d-c\} = N(-c).\]
Therefore $W$ is orthogonal to $\col\pmat{A_1 & A_2}$. Finally, since $X$ is $4$-regular, by Corollary \ref{twin-supp}, $x_a(W)$ and $x_b(W)$ are strongly-cospectral with
\[\Lambda^+_{a(W), b(W)} = \{\pm 1/\sqrt{2} \}= \{\cos(\pi/4),  \cos(3\pi/4)\},\]
and 
\[\Lambda^-_{a(W), b(W)}=\{0\}=\{\cos(2\pi/4)\}.\]
It follows from Theorem \ref{char-twins} (i) that pointwise perfect $W$-transfer occurs from $a$ to $b$ at time $4$.
\qed		

A \textsl{double cone} over $X$, denoted $\comp{K_2}+X$, is a graph obtained from $X$ by adding two vertices, called the \textsl{conical vertices}, and joining them to all vertices in $X$. The following result shows a double cone over the disjoint union of $k$ cycles, whose lengths are all divisible by $4$, can transfer a $2k$-dimensional subspace between the two marked conical vertices.

\begin{theorem}
Let $X$ be the double cone
	\[X = \comp{K_2} + ( C_{4m_1} \cup  C_{4m_2} \cup \cdots \cup  C_{4m_k})\]
	where $\seq{m}{1}{2}{k}$ are some positive integers.
	Let $a$ and $b$ be the conical vertices. 	Let $W$ be the subspace of functions on $V(X)\setminus\{a,b\}$ spanned by those that take alternating values $1,0,-1,0$ along some cycles $C_{4m_j}$ and vanish on all other cycles. Let $C_a=C_b$ be the reflection about any subspace containing $W$. There is pointwise perfect $W$-transfer from $a$ to $b$ at time $t=4$.
\end{theorem}
\proof
Any function that takes alternating values $1,0,-1,0$ along $C_{4m_j}$ lies in the kernel of $A(C_{4m_j})$. Hence, $W$ is orthogonal to the column space of $A(X\setminus\{a,b\})$. Since all neighbors of $a$ have degree $4$, by Corollary \ref{twin-supp}, $x_a(W)$ and $x_b(W)$ are strongly-cospectral with
\[\Lambda^+_{a(W), b(W)} = \{\pm 1/\sqrt{2} \}= \{\cos(\pi/4),  \cos(3\pi/4)\},\]
and 
\[\Lambda^-_{a(W), b(W)}=\{0\}=\{\cos(2\pi/4)\}.\]
It follows from Theorem \ref{char-twins} (i) that pointwise perfect $W$-transfer occurs from $a$ to $b$ at time $4$.
\qed

If we relax to pointwise pretty good subspace state transfer, we obtain more examples from double cones. To start, we cite a useful result about geodetic angles. A real number is a \textsl{pure geodetic angle} if any one of its six squared trigonometric functions is either rational or infinite. 

\begin{theorem}\cite[Ch 11]{Bergen2009}\label{thm:tan_rat_pi} % 
	If $q\in \rats$ and $q\pi$ is a pure geodetic angle, then 
	\[\tan(q\pi) \in \left\{0, \pm{\sqrt{3}}, \pm \frac{1}{\sqrt{3}}, \pm 1\right\}.\]
\end{theorem}

\begin{theorem}
	Let $Z$ be a $k$-regular graph with singular adjacency matrix. Let $W$ be a subspace of $\ker(A(Z))$. Let $X=\comp{K_2}+Z$ with conical vertices $a$ and $b$. Let $C_a = C_b$ be the reflection about any subspace containing $W$. If $k\notin\{0,2,6\}$, then pointwise pretty good $W$-transfer occurs from $a$ to $b$ in integer steps.
\end{theorem}
\proof
By Corollary \ref{twin-supp} we see that $x_a(W)$ and $x_b(W)$ are strongly-cospectral, with
	\[ \Lambda_{a(W), b(W)}^+ =\left\{\pm \sqrt{\frac{2}{k+2}}\right\}, \quad \Lambda_{a(W), b(W)}^- = \{0\}.\]
It suffices to show Theorem \ref{basis-free-pgst} (ii) holds for $\gamma=-1$.  Let $\ell_0$, $\ell_+$ and $\ell_-$ be integers such that 
\[\ell_0 \arccos 0 + \ell_+ \arccos\sqrt{\frac{2}{k+2}} + \ell_- \arccos \left(-\sqrt{\frac{2}{k+2}}\right)  \equiv 0\pmod{2\pi}.\]
Then
\[\left(\frac{\ell_0}{2}+\ell_-\right) \pi + (\ell_+ - \ell_-) \arccos\sqrt{\frac{2}{k+2}} \equiv 0\pmod{2\pi}.\]
If $\ell_+\ne \ell_-$, then 
\[\arccos\sqrt{\frac{2}{k+2}} = q\pi \]
for some $q\in\rats$, and $q\pi$ would be pure geodetic. However, Theorem \ref{thm:tan_rat_pi} implies this is impossible for $k\notin\{0,2,6\}$. Therefore $\ell_+=\ell_-$, and we have
\[\sum_{\lambda\in \Lambda_{a(W), b(W)}^+} \ell_{\lambda} = \ell_++\ell_-\equiv 0\pmod{2}. \tag*{\sqr53}\]

We now consider another special case where $N(a)$ and $N(b)$ are disjoint. As before, assume $a$ and $b$ are not adjacent. Then
\begin{center}
	\begin{tikzpicture}
		\draw [thick](0.25,0) -- (0.25,5.5);
		\draw [thick](0.25,0) -- (0.5,0);
		\draw [thick](0.25,5.5) -- (0.5,5.5);
		\draw [thick](7,0) -- (7,5.5);
		\draw [thick](6.75,0) -- (7,0);
		\draw [thick](6.75,5.5) -- (7,5.5);
		\draw [thick] (3,0) rectangle (6.5,2.75);
		\node [scale=1] at (4.75,1.375) {$A(X\setminus\{a,b\})$};

		\draw [thick,color=blue!60] (0.5,1.75) rectangle (1.5,2.75);
		\node [scale=1,color=blue] at (1,2.25) {$\frac{C_a+I}{2}$};
		\draw [thick,color=red!60] (1.75,0.75) rectangle (2.75,1.75);
		\node [scale=1,color=red] at (2.25,1.25) {$\frac{C_b+I}{2}$};
		
		\draw [thick,color=blue!60] (3,4.25) rectangle (4,5.25);
		\node [scale=1,color=blue] at (3.5,4.725) {$\frac{C_a+I}{2}$};
		\draw [thick,color=red!60] (4,3) rectangle (5,4);
		\node [scale=1,color=red] at (4.5,3.5) {$\frac{C_b+I}{2}$};

		\node [scale=1] at (2.25,0.25) {$0$};
		\node [scale=1] at (1,0.25) {$0$};
		
		\node [scale=1] at (5.75,4.725) {$0$};
		\node [scale=1] at (5.75,3.5) {$0$};
		
		\node [scale=1] at (1,4.75) {$0$};
		\node [scale=1] at (2.25,4.75) {$0$};
		\node [scale=1] at (1,3.5) {$0$};
		\node [scale=1] at (2.25,3.5) {$0$};
		
		\node [scale=0.7,color=blue] at (1,5.75) {$\cl(a)$};
		\node [scale=0.7,color=red] at (2.25,5.75) {$\cl(b)$};
		
		\node [scale=0.7,color=blue] at (3.5,5.75) {$N(a)$};
		\node [scale=0.7,color=red] at (4.5,5.75) {$N(b)$};
	
\node [scale=1] at (-1.2,5.25) {$G=\pmat{I & \\ & \Delta}^{-\frac{1}{2}}$};
\node [scale=1] at (8,5.25) {$\pmat{I & \\ & \Delta}^{-\frac{1}{2}}$};
	\end{tikzpicture}
\end{center}

		We use this observation to construct an infinite family of graphs with pointwise subspace state transfer between vertices at arbitrary distance. Let $GP(k, n)$ denote the graph obtained from the disjoint union of $k$ paths $P_n$ by identifying their left endpoints and right endpoints respectively, as shown in Figure \ref{GP}.
		
		\begin{figure}[h]
			\centering
		        \begin{tikzpicture}
			[scale=1, every node/.style={circle,very thick,draw}]
			\tikzset{lab/.style={draw=none,  inner sep = 4pt}}
			
			\draw[fill=gray!45](0,0) ellipse [x radius=3.75, y radius=0.75];
			\draw[fill=gray!25, fill opacity=0.3](-1.5,0) ellipse [x radius=0.75, y radius=3];
			
			\node[fill=red] (0) at (3,0) {};
			\node[fill=blue] (1) at (-3,0) {};
			\node[] (2) at (1.5,0) {};
			\node[] (3) at (1.5,2) {};
			\node[] (4) at (1.5,-2) {};
			\node[] (5) at (0,0) {};
			\node[] (6) at (0,2) {};
			\node[] (7) at (0,-2) {};
			\node[] (8) at (-1.5,0) {};
			\node[] (9) at (-1.5,2) {};
			\node[] (10) at (-1.5,-2) {};
			
			\node[draw=none, fill=none] at (-0.75,0) {\textbf{\Large$\dots$}};
			\node[draw=none, fill=none,] at (4.25,0) {\textbf{\Large$n$}};
			
			\node[draw=none, fill=none] at (-1.5,1) {\Large$\vdots$};
			\node[draw=none, fill=none,] at (-1.5,3.5) {\textbf{\Large$k$}};

			\node[draw=none, fill=none,text=blue] at (-3.5,0) {\Large$a$};
			\node[draw=none, fill=none,text=red] at (3.5,0) {\Large$b$};

			\draw [thick](1) -- (-1.6,1);
			\draw [thick](-1.4,1) -- (1.5,1);
			\draw [thick](1.5,1) -- (0);
			
			\foreach \a/\b in {1/8,1/9,1/10,9/6,10/7,5/2,6/3,7/4,2/0,3/0,4/0}
			\draw[thick] (\a) to (\b);	
		\end{tikzpicture}
		\caption{$GP(k,n)$}
		\label{GP}
		\end{figure}
		
		\begin{theorem}
			Let $X=GP(k,n)$ with $a$ and $b$ being the vertices of degree $k$. Let $U=RC$ be the transition matrix of a Grover coined walk on $X$ with marked vertices $a$ and $b$ and any reflection coin $C_a=C_b$. Let $W$ be any subspace of $\col(C_a+I)$. There is pointwise perfect $W$-transfer from $a$ to $b$ at time $n-1$.
		\end{theorem}
		\proof
	Let $Q$ be the projecion onto $\col(C_a+I)$. Since $X\backslash\{a,b\}$ is a disjoint union of $k$ paths $P_{n-2}$, we can rearrange $G$ as
\[
G =
\pmat{
	0 & \frac{1}{\sqrt{2}}Q & 0 & 0 & \cdots & 0 & 0 \\
	\frac{1}{\sqrt{2}}Q & 0 & \frac{1}{2}I& 0 & \cdots & 0 & 0 \\
	0 & \frac{1}{2}I & 0 & \frac{1}{2}I & \cdots & 0 & 0 \\
	0 & 0 & \frac{1}{2}I & 0 & \ddots & \vdots & \vdots \\
	\vdots & \vdots & \vdots & \ddots & \ddots & \frac{1}{2}I & 0 \\
	0 & 0 & 0 & \cdots & \frac{1}{2} I & 0 & \frac{1}{\sqrt{2}} Q \\
	0 & 0 & 0 & \cdots & 0 & \frac{1}{\sqrt{2}} Q& 0
}.
\]
Note that $G$ acts on $\cx^n\otimes W$ as the normalized adjacency matrix of $P_n$, which has eigenvalues and eigenvectors
\[
	\lambda = \cos \frac{\ell \pi}{n-1},\qquad z = \pmat{ \frac{1}{\sqrt{2}} \\ \cos \frac{\ell \pi}{n-1}\\\cos \frac{2\ell \pi}{n-1}\\\vdots \\\cos \frac{(n-2)\ell \pi}{n-1}\\ \frac{(-1)^{\ell}}{\sqrt{2}}},\qquad \ell=0,1,\cdots,n-1.
\]
Now let $z$ be any eigenvector for $G$ with eigenvalue $\lambda$. We can decompose $z$ into $z = x+y$
where $x\in \cx^n\otimes W$ and $y\in \cx^n \otimes W^{\perp}$. Thus either 
\[\projonto{W}{z[\cl(a)]} = \projonto{W}{z[\cl(b)]} =0,\]
or
\[\lambda  = \cos \frac{\ell \pi}{n-1},\quad \projonto{W}{z[\cl(a)]} = (-1)^{\ell} \projonto{W}{z[\cl(b)]}\ne 0.\]
It follows from Theorem \ref{basis-free} that pointwise perfect $W$-transfer occurs from $a$ to $b$ at time $n-1$.
\qed

\section{Discussion}		
We have developed some theory for perfect and pretty good subspace state transfer on simple unweighted graphs. Our characterizations connect these transport properties of walks to spectral properties of certain Hermitian weighted digraphs. An interesting direction is to explore spectral properties of the $(C_a, C_b)$-blow-ups, and use them to classify graph structures that admit subspace state transfer relative to specific reflection coins such as negative Grover coins and Hadamard coins. We plan on investigating these questions in the future.

\section*{Acknowledgement}
This material is based upon work supported by the National Science Foundation under Grant No. 2348399.

		\bibliographystyle{amsplain}
		\bibliography{qw}
		
	\end{document}